\documentclass[11pt]{article}

\usepackage[alldates=year, block=space, style=apa]{biblatex}
\usepackage[margin=1in]{geometry}
\usepackage[colorlinks, allcolors=blue]{hyperref}

\usepackage{amsmath, amssymb, amsthm, amsfonts, authblk, color, enumerate, framed, graphicx, mathrsfs, mathtools, mdframed, subcaption, tabularx, theoremref, titlesec, verbatim, xfrac, soul, mathtools, wrapfig, tikz, float}

\graphicspath {{figures/}}

\newtheorem*{theorem*}{Theorem}
\newtheorem{theorem}{Theorem}
\newtheorem*{corollary*}{Corollary}

\newtheorem*{lemma*}{Lemma}
\newtheorem{lemma}{Lemma}
\newtheorem*{conjecture*}{Conjecture}

\newtheorem{proposition}{Proposition}
\newtheorem*{assumption*}{Assumption}
\newtheorem{assumption}{Assumption}
\newtheorem*{definition*}{Definition}
\newtheorem{definition}{Definition}

\usepackage[colorinlistoftodos]{todonotes}
\usepackage[normalem]{ulem}

\newcommand{\argmax}{\operatornamewithlimits{argmax}}
\newcommand{\argmin}{\operatornamewithlimits{argmin}}
\newcommand{\defeq}{\stackrel{\text{def}}{=}}
\newcommand{\abs}[1]{\left\lvert#1\right\rvert}

\title{
  Strategies under a pickoff limit in baseball: A zero-sum sequential game with multilevel models
}

\author[1]{Scott Powers}
\author[2]{Sivaramakrishnan Ramani}
\author[3]{Jacob Hahn}
\author[2]{Andrew J. Schaefer}
\affil[1]{\small Department of Sport Management, Rice University}
\affil[2]{\small Department of Computational Applied Mathematics and Operations Research, Rice University}
\affil[3]{\small Department of Statistics, Rice University}

\date{}

\begin{document}

\maketitle

\begin{abstract}
  Recently, Major League Baseball (MLB) limited pitchers to three pickoff attempts, creating a cat-and-mouse game between pitcher and runner. Each failed attempt adds pressure on the pitcher to avoid using another, and the runner can intensify this pressure by extending their leadoff toward the next base. We model this dynamic as a two-player zero-sum sequential game in which the runner first chooses a lead distance, and then the pitcher chooses whether to attempt a pickoff. We establish optimality characterizations for the game and present variants of value iteration and policy iteration to solve the game. Using high-resolution lead distance data from optical player tracking, we estimate generalized linear mixed-effects models for pickoff and stolen base outcome probabilities given lead distance, context, and player skill. We compute the game-theoretic equilibria under the two-player model, as well as the optimal runner policy under a simplified one-player Markov decision process (MDP) model. In the one-player setting, our results establish an actionable rule of thumb: the Two-Foot Rule, which recommends that a runner increase their lead by two feet after each pickoff attempt.
\end{abstract}

 \section{Introduction}
  \label{sec:introduction}
    On May 17, 2025, the Chicago Cubs hosted the Chicago White Sox in the second game of their Crosstown Series. In the bottom of the first inning, White Sox pitcher Sean Burke faced Cubs batter Michael Busch with two outs and Cubs runner Kyle Tucker on first base. Before the 1-1 pitch, Burke pivoted and threw to first base behind Tucker, who retreated safely before a tag. Tucker had 10 stolen bases in 10 attempts on the young season, the most of any base stealer who had not been caught \parencite{fangraphs_major_2025}. On the White Sox broadcast, color commentator Steve Stone explained to the audience, ``The only problem is now there's two throws over there, and Tucker's gonna get a much bigger lead.'' On the next pitch, Tucker stole his 11$^\text{th}$ base of the season, the first Burke had allowed. ``I don't think anybody in the room could have thrown him out on that play,'' said Stone.

    For background, the stolen base is an exciting play in which a runner takes off on the pitch, attempting to reach the next base before the catcher can throw the ball to a fielder at the base to tag the runner out. Whether attempting a stolen base or not, the runner takes a {\it lead} from their base, commonly positioning themselves 8--12 feet away from their base, toward the next base. The longer the lead, the shorter the distance needed to run to reach the next base. The pitcher's primary deterrent against long leads is a {\it pickoff attempt}, in which they throw the ball to a fielder at the runner's base, attempting a tag before the runner returns safely to the base.

    Major League Baseball (MLB) stolen base totals surged in the 1980s and 1990s but then waned in the 2000s and 2010s as teams cut back on the play, which generally requires a high success rate in order to have a positive impact on run scoring \parencite{tango_book_2007}. Before 2023, pickoff attempts could delay games because unsuccessful pickoff attempts do not advance the state of the game. Acknowledging this and other trends, MLB introduced several rule changes in 2023, with the goal of making the game more exciting \parencite{castrovince_pitch_2023}. With these changes, a pitch clock limits the time between pitches, and there are restrictions on infielder positioning, curtailing defensive shifts. Two changes are particularly relevant to stolen base attempts: Bases are now wider (reducing the distance between first and second base by 4.5 inches), and there are now limitations on pickoff attempts, or more formally, {\it disengagements}.

    Under the new rules, pitchers are limited to three disengagements. When attempting a pickoff, instead of throwing a pitch, the pitcher disengages from the pitching rubber and throws the ball to the runner's base. If the pitcher disengages from the pitching rubber three times without successfully picking off the runner, then all runners automatically advance to the next base. Disengagements also include any time the pitcher steps off the rubber without attempting a pickoff, which we disregard in the present work. The disengagement count resets if a runner advances or if a plate appearance ends.

    This disengagement limit creates an interesting cat-and-mouse game between the pitcher and the runner. Each unsuccessful pickoff attempt disincentivizes the pitcher from attempting another pickoff, which incentivizes the runner to increase their lead distance. As indicated by Stone before the Tucker stolen base against the White Sox, players and coaches understand this intuitively. But how much should the runner extend their lead, and when should the pitcher attempt a pickoff?

    \subsection{Contributions}

\begin{figure}
      \centering
      \includegraphics[width = 0.8\linewidth]{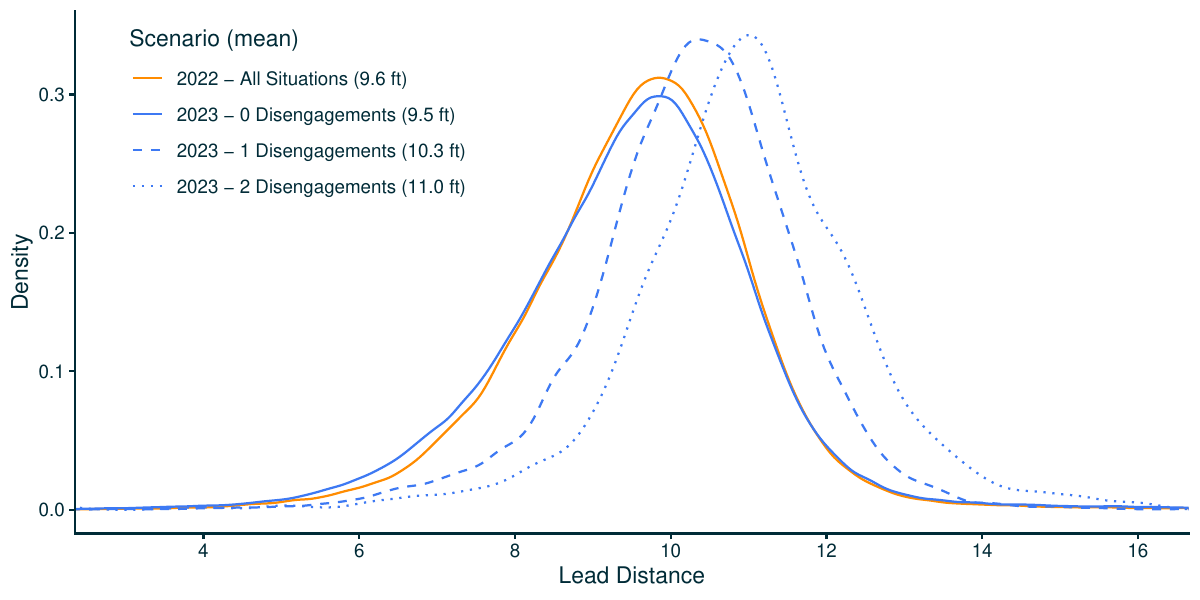}
      \caption{
        \it The league-wide distribution of lead distance from first base when second base is unoccupied, comparing 2022 with 2023. The 2023 data are split by prior disengagements because, beginning in 2023, pitchers were limited to three disengagements, influencing pitcher and runner behavior. In 2022, there was no limit on the number of disengagements.
      }
      \label{fig:leads-overall}
    \end{figure}

    To address the problem of how to select a lead distance under a pickoff limit in baseball, our work leverages newly observable continuous spatial tracking data from MLB. Derived from optical player tracking, this granular dataset comprises hundreds of thousands of high-fidelity observations of exact lead distances taken by runners, measured via center-of-mass projection. The data covers all MLB regular season games in 2022 and 2023, spanning both sides of the pickoff limit rule change. Figure \ref{fig:leads-overall} illustrates that in 2023, before the first disengagement, the distribution of lead distances was similar to 2022. Each disengagement corresponded to a successive 0.7-foot increase in mean lead distance.

    We introduce a formal model of the baseball inning as a two-player zero-sum sequential game. The model captures the vital aspects of the runner-pitcher interaction while retaining analytical and computational tractability. Our model is closely related to the class of stochastic/dynamic Stackelberg games \parencite{goktas_zero-sum_2022,li_review_2017,lopez_stationary_2023,vasal_sequential_2022,vorobeychik_computing_2012}, which combines the elements of Stackelberg and stochastic games. The main difference between our infinite-horizon model and the stochastic Stackelberg games literature is that we do not discount the rewards. The rationale for this modeling choice is that our optimization objective---expected runs in the inning---by definition includes no discount. We establish a maximin Bellman equation for our sequential game and propose a variant of value iteration and policy iteration to determine the optimal strategies for runner and pitcher.

    Rather than including player identities in the state space, which would cause severe data sparsity issues, we account rigorously for individual player skills by employing a two-step modeling procedure to infer the transition probabilities. First, using the law of total probability, we decompose the global transition function into granular, component runner outcomes. Second, we estimate probability models for component runner outcomes using generalized linear mixed-effects models (GLMMs). In addition to conditioning on game context and player actions, the multilevel models incorporate both observable player physical traits via fixed effects and unobserved player latent skills via random effects. We multiply the component probability estimates with pooled conditional empirical frequencies to assemble transition probabilities for downstream policy evaluation.

    Through a comparative analysis of observed player behavior and the game-theoretic equilibrium, we show that both pitcher and runner actions are very conservative relative to optimal policies. Pitchers attempt fewer pickoffs than optimal while runners take shorter leads than optimal, and moving to a rational equilibrium would benefit the offense. Finally, we consider a reduction of the sequential game to a (one-player) Markov decision process (MDP) and compress the high-dimensional optimal policy into a highly actionable heuristic for a runner seeking to exploit pitcher behavior: The Two-Foot Rule, which recommends that the runner increase their lead by two feet after each pickoff attempt.

    While our application is specific to the problem of selecting a runner's lead distance in baseball, we develop a template that is transferable more broadly. Sports officials often consider rule changes that impact the structure of competition, such as the kickoff in the National Football League \parencite{yousuf_nfls_2026} or the draft lottery in the National Basketball Association \parencite{vorkunov_nbas_2026}. Our framework is useful both to teams seeking to maximize their success under changing rules and to leagues seeking to understand the potential consequences of rule changes before they are implemented.

    \subsection{Related work}

    \textcite{downey_pick_2015} analyzed pickoff throws and stolen base attempts as a mixed strategy game without using lead distance data. They modeled the two-player game as a simultaneous decision in which the pitcher chooses whether to attempt a pickoff and the runner chooses whether to attempt a stolen base, showing that this theoretical game has a mixed-strategy Nash equilibrium. Under this equilibrium, left-handed pitchers (typically possessing better pickoff moves than their right-handed counterparts) attempt fewer pickoffs and see fewer stolen bases attempted against them, a result validated by empirical data. In a later analysis, \textcite{downey_pressure_2019} argued that pitchers can effectively implement a mixed strategy in low-pressure situations but struggle to randomize their decisions as pressure increases, showing a higher auto-correlation for sequential pickoff decisions in high-pressure situations. Both of these analyses covered a game in which pickoff attempts were unlimited.

    Several studies in the sports economics literature have examined the extent to which professional athletes play minimax strategy. \textcite{palacios-huerta_professionals_2003} found evidence consistent with minimax equilibrium behavior in the mixed strategy game of the soccer penalty kick, in which the kicker chooses to shoot the ball left or right and the goalkeeper simultaneously chooses to dive left or right for the save. \textcite{palacios-huerta_experientia_2008} followed this with a laboratory experiment featuring professional soccer players and college students simulating the penalty kick decision using cards, showing that the professionals played minimax but the students did not. On the other hand, \textcite{kovash_professionals_2009} found evidence inconsistent with minimax equilibrium behavior in the National Football League (teams called more run plays and fewer pass plays than theoretically optimal) and in MLB (pitchers threw more fastballs than theoretically optimal).

    The history of quantitative analysis in baseball has deep roots in operations research (e.g., \cite{lindsey_investigation_1963}; \cite{freeze_analysis_1974}; \cite{pankin_evaluating_1978}; \cite{russell_devising_1994}; \cite{bukiet_markov_1997}). The Markov process is particularly well suited for baseball, where the game's highly structured progression lends it well to an intuitive state space model \parencite{sokol_intuitive_2004}. \textcite{hirotsu_markov_2003} used a Markov model to inform substitution decision-making (specifically, pinch hitting) in baseball. The Markov process is also popular across sport analytics. \textcite{hirotsu_using_2002} used a four-state Markov model for soccer to inform substitution decision-making based on how substitutions impact transition probabilities between states. Recently, \textcite{chan_points_2021} used a Markov model to evaluate team performance in football. More specifically, MDP has become a particularly popular approach in recent years for analyzing decision-making in sports, with applications in tennis \parencite{nadimpalli_when_2013}, basketball \parencite{sandholtz_markov_2020}, darts \parencite{baird_optimising_2020}, and football \parencite{biro_reinforcement_2022}. In one particularly relevant baseball application, \textcite{hirotsu_using_2019} used an MDP to model sacrifice bunt decisions on the basis of win probability, finding that the sacrifice bunt was more valuable than previously thought. Sequential games are emerging as a preferred approach to analyzing competition in sports. \textcite{haugh2022play} analyzed competition in darts using a two-player zero-sum stochastic game with undiscounted rewards. They exploited the finiteness and state-independence of the players' action sets to reformulate the alternate-move darts game as a zero-sum (simultaneous-move) stochastic game and then analyzed it using the theory developed in the literature for stochastic games. Applying the reformulation tricks of \textcite{haugh2022play} for our baseball problem would render it as a stochastic game with an uncountable state space, which would be analytically and computationally intractable.

    Earlier work proposed probability models for stolen bases in contexts different than the one considered in this paper. \textcite{loughin_assessing_2008} estimated linear mixed-effects models for stolen base attempts and successes, using random effects for pitcher and catcher (not the runner, notably) and fixed effects for balls, strikes, outs, score, venue and pitcher hand. They showed that the spread of pitcher effects is slightly greater than the spread of catcher effects. Unlike our paper, that analysis did not use data on runner sprint speed and catcher arm strength. More recently, \textcite{stanley_modeling_2023} modeled stolen base success probability using both logistic regression and a random forest, including runner sprint speed and catcher arm strength as features but excluding player identities. Furthermore, that analysis was limited to data prior to the 2023 MLB rule changes that incentivized base stealing. We emphasize that the probability models in Section \ref{sec:transition-probability-model} extend these prior works by combining fixed effects for runner sprint speed and catcher arm strength with player random effects, as well as incorporating lead distance and modeling pickoff attempts and successes.
  
  \subsection{Data}

    The new data which make it possible to address the lead distance problem are play-by-play lead distances measured via optical player tracking. The runner's lead distance is defined as the distance between the runner's center of mass (projected onto the line between the first base and second base) and the nearest edge of first base. This distance is measured at the time of the pitcher's first movement, which could be a pitch or a pickoff attempt. We obtained runner lead distances from MLB Advanced Media for every play in the 2022 and 2023 MLB regular seasons, covering the season before and the season after the rule change introducing the pickoff limit. When the pitcher throws a pitch, the runner extends their lead by several feet from the time of the pitcher's first movement to the time of pitch release (known as taking a ``secondary lead''). However, the data include only the lead distance at the time of the pitcher's first movement, not the secondary lead.

    In addition to the proprietary lead distance data, we downloaded publicly available contextual data from the MLB Stats API using the R package sabRmetrics \parencite{powers_sabrmetrics_2025}. The dataset comprises 722,899 plays from the 2022 MLB regular season and 726,868 plays from the 2023 MLB regular season, where each play is either a pitch or a pickoff attempt (2.0\% pickoff attempts in 2022, 1.2\% in 2023). For each play, the data include bases occupied, ball-strike count, disengagements, and outs, as well as identities of the pitcher, catcher, and runners involved. On each pickoff attempt, we observe whether it is successful. On each pitch, we observe whether the runner attempts a stolen base; and if so, we observe whether they are successful. On each pitch terminating a plate appearance, we observe the plate appearance outcome, such as strikeout, groundout, single, or home run.

    Finally, we obtained rich, publicly available data on catcher and runner skills measured via optical player tracking \parencite{baseball_savant_statcast_2024}. For each catcher-season, we observe {\it arm strength} (miles per hour), which Baseball Savant defines as the average of the top 10 percent (to capture peak performance) of throw speeds when the catcher is attempting to catch a base stealer. For each runner-season, we observe {\it sprint speed} (feet per second), which Baseball Savant defines as the average of the top 67 percent (to exclude low-effort plays) of peak running speeds on qualified runs: (1) home-to-first runs on weakly hit batted balls and (2) runs of two or more bases on balls hit into play (excluding runs from second base on extra-base hits).

  \subsection{Outline}

    The rest of this paper is organized as follows. In Section \ref{sec:sg-model}, we propose a theoretical model for the dynamic between runner and pitcher under the pickoff limit rule as a two-player zero-sum sequential game. In Section \ref{sec:transition-probability-model}, we propose a hierarchical model for estimating transition probabilities in the sequential game. In Section \ref{sec:results}, we present the results of estimating the transition probability model and solving the sequential game using data from the 2023 MLB regular season. In Section \ref{sec:discussion}, we discuss actionable insights from the results, including the Two-Foot Rule, which is based on reducing the two-player sequential game to a one-player MDP. Finally, in Section \ref{sec:conclusion}, we conclude with a summary of the key results and their significance.

  \section{Two-player sequential game model}
  \label{sec:sg-model}

    We begin with a model of a baseball inning that captures the vital aspects of the pitcher-runner game to enable the generation of actionable insights for the runner and the pitcher. We take runs as the objective so that the runner aims to maximize expected runs scored until end of inning and the pitcher aims to minimize the same. Although game win probability would be a more pure objective, the choice of inning run expectation leads to more actionable insights, as we discuss in Section \ref{sec:discussion}. Our model assumes a representative pitcher and a representative runner for the whole inning, rather than capturing full lineups and substitutions of players with varying skills. Without this assumption, we would end up with a nonstationary (infinite-horizon) two-player sequential game, which appears to be significantly more challenging to analyze when compared with the stationary version.

    \subsection{Components}
    \label{sec:sequential_gama_components} 
    The four components of our sequential game are as follows.

    \noindent{\bf 1. State space.} We use $s=(b,c,d,o) \in {\cal X}$ to capture a play, where     
    \begin{itemize}
      \item $b \in \{0, 1\}^3$ are the bases occupied by runners;
      \item $c \in \{0, 1, 2, 3\} \times \{0, 1, 2\}$ is the ball-strike count;
      \item $d \in \{0, 1, 2\}$ denotes the number of disengagements already used; and
      \item $o \in \{0, 1, 2\}$ denotes the number of outs.
    \end{itemize}
    The number of elements in ${\cal X}$ equals $2^3 \times (4 \times 3) \times 3 \times 3 = 864$. We use $\widetilde{\cal X}$ to denote the set of four penultimate states and $\Delta$ to denote the terminal end-of-inning state. The purpose of the penultimate states is to keep track of the number of runs that are scored on the transition to the terminal state. A penultimate state $s \in \widetilde{\cal X} \equiv \{0, 1, 2, 3\}$ represents the number of runs scored on the final play of the inning. We use $\mathcal{S} \equiv \mathcal{X} \cup \widetilde{\cal X} \cup \{\Delta\}$ to denote the set of all states and $S_t$ to denote the state at play $t$ in a (fixed) inning. Traditionally, the state is defined by the bases occupied and the number of outs (24 non-terminal states) or by those and the ball-strike count (288 non-terminal states) \parencite{bukiet_markov_1997}. Our construction of the state extends them by further including the number of disengagements as a component of the state, thereby allowing us to analyze the game under the new pickoff limit rule.

   \noindent{\bf 2. Action space.} We limit the agency of runners and pitchers to states in which only first base is occupied (i.e., the component $b$ in the state variable equals $(1,0,0)$). Although this restriction limits the scope, these states covered roughly $71\%$ of stolen base attempts in 2023. We emphasize that this  excludes the scenario with $b = (1, 0, 1)$ in which the first-base runner will often attempt a stolen base to lure a throw from the catcher, giving the third-base runner a chance to steal home and establishing a different dynamic \parencite{miller_first-and-third_2024}. We leave this for future research.
    
    For states in which only first base is occupied, the first-base runner's action is the lead distance from first base, chosen in feet. The pitcher has only two actions: attempt a pickoff or throw a pitch. Formally, the action space of the first-base runner, $A_{\text{R}}$, and the action space of the pitcher, $A_{\text{P}}$, are as follows.
    
    \begin{align}
    \label{eqn:runner_action_space}
        A_\text{R}(s) \defeq \begin{cases}
            {\cal L} \text{ if } s = (b,c,d,o) \text{ with } b = (1,0,0)\\
            \{\delta\} \text{ if } s = (b,c,d,o) \text{ with } b \neq (1,0,0), \text{ or } s \in \widetilde{\cal X} \cup \{\Delta\}
        \end{cases} \forall s \in {\cal S},
    \end{align}
        
    where ${\cal L} \defeq [0,20]$ is the set of allowable lead distances for the first-base runner, and
    
    \begin{align}
      \label{eqn:pitcher_action_space}
      \begin{split}
        A_{\text{P}}(s,a_\text{R}) \defeq&
        \begin{cases}
          {\cal P} \text{ if } s = (b,c,d,o) \text{ with } b = (1,0,0) \\
          \{\delta\} \text{ if } s = (b,c,d,o) \text{ with } b \neq (1,0,0), \text{ or } s \in \widetilde{\cal X} {\cup} \{\Delta\},
        \end{cases}\\
        & \qquad \forall s \in {\cal S}, a_\text{R} \in A_\text{R}(s),
      \end{split}
    \end{align} 
    where ${\cal P} \defeq \{0,1\}$, with $p = 1$ corresponding to a pickoff attempt and $0$ otherwise.  Here, $\delta$ is a dummy element that captures the lack of agency. We note that the state of our sequential game progresses even when no actions are played by both players. We will model this through appropriate transition probabilities.  

    \noindent{\bf 3. Transition probabilities.} The uncertainty in our sequential game is modeled through transition probabilities, i.e., the conditional probability of transitioning to a future state given the sequence of previously visited states and the actions played in those states. Our sequential game assumes that the transition probabilities satisfy the Markov property.  Following the convention in the literature, we denote the transition probability as $p(\cdot \mid s,a_\text{R},a_\text{P}) \ \forall s \in {\cal S}, a_\text{R} \in A_\text{R}(s), a_\text{P} \in A_\text{P}(s,a_\text{R})$. Section \ref{sec:transition-probability-model} details our specific construction of transition probabilities, which uses generalized linear mixed-effects models to estimate runner outcome (e.g., unsuccessful pickoff attempt or successful stolen base attempt) probabilities which depend on lead distance and the pitcher's pickoff decision, as well as player-specific factors. These, along with appropriate conditional empirical frequencies, will be used to construct the transition probabilities.

    \noindent{\bf 4. Reward function.} The final component of our sequential game is the reward function, which models the payoff received by the players. Since we use a  zero-sum sequential game, the payoff obtained by one player comes at the expense of the other player. We choose the number of runs scored as the reward function. This design choice is more suitable for the early and middle innings of a game and less suitable for the later innings, particularly the bottom of the ninth or extra innings, when the offense knows how many runs they need to score to win the game. For any state $s \in {\cal S}$, the reward function captures the number of runs scored by the batting team in transitioning from state $s$ to $s'$ given the actions of the first-base runner and the pitcher. Formally,
    \begin{align}
      \label{eqn:runner_pitcher_reward}
        r(s'|s,a_\text{R},a_\text{P}) &\defeq r(s'|s) 
        \, \forall s,s' \in {\cal S},\, a_\text{R} \in A_\text{R}(s),\,  a_\text{P} \in A_\text{P}(s,a_\text{R}),
    \end{align}
    where
    \begin{align}
      \label{eqn:common_reward}
      r(s'|s) \defeq  \begin{cases}
          (g(b) + o) - (g(b') + o') +\mathbb{I}\{c' = (0, 0),\, d' = 0\} \text{ if } s,\, s' \in \mathcal{X}\\
          s'  \text{ if } s \in \mathcal{X},\, s' \in \widetilde{\cal X}\\
          0  \text{ if } s' = \Delta,
      \end{cases} 
    \end{align}
    where $g: \{0, 1\}^3 \rightarrow \{0, 1, 2, 3\}$ counts the number of runners on base. Note that our reward function in \eqref{eqn:runner_pitcher_reward} depends only on the current state and the next state, and not on the actions of the first-base runner or the pitcher. The formula in \eqref{eqn:common_reward} infers the number of runs scored by counting the offensive players captured by the states before and after the transition. If the transition corresponds to the end of a plate appearance (i.e. $c' = (0, 0),\, d' = 0$), then we account for one more offensive player (the batter) in the next state. Since the number of runs scored on the last play cannot be uniquely determined from the current state, we use the penultimate states to track the number of runs scored on the last play. We note that our transition probabilities will be constructed such that these penultimate states transition to the terminal end-of-inning state with probability $1$. Since the reward function \eqref{eqn:runner_pitcher_reward} depends only on the current and the next state, it appears to be independent of the actions of the players. However, the next state to which our sequential game transitions is (probabilistically) determined by the actions of both the players through the transition probabilities, thereby making the rewards depend indirectly on the actions of both players.

    \subsection{Formulating and solving the sequential game}
    \label{sec:sequential_game_Setup}
    Our sequential game models the progression of an inning as follows. Let $s \in {\cal S}$ denote the state of the sequential game from which a play begins. Upon observing this state, the first-base runner picks an action $a_\text{R} \in A_\text{R}(s)$. The pitcher observes the state $s$ and the first-base runner's action $a_\text{R}$, and then picks an action $a_\text{P} \in A_\text{P}(s,a_\text{R})$. After playing $a_\text{R}$ and $a_\text{P}$, the sequential game then transitions into state $s'$ with probability $p(s'|s,a_\text{R},a_\text{P})$ and the batting team collects a reward $r(s'|s,a_\text{R},a_\text{P})$ in the zero-sum game. A play then begins from $s'$, and the process repeats. The goal of the batting team is to maximize (through the actions of the first-base runner) the expected total rewards over the plays in the inning while the goal of the fielding team is to minimize (through the actions of the pitcher) the expected total rewards obtained by the batting team.

    Intuitively, any sound strategy for the first-base runner and the pitcher depends on the state of the sequential game (for the pitcher, it also depends on the actions played by the first-base runner). This is formally captured using the concept of policy. A (deterministic, stationary) policy for the first-base runner is a decision-rule that prescribes the action to be played in every possible state. Similarly, a (deterministic, stationary) policy for the pitcher is a decision-rule that prescribes the action to be played for every possible state and the corresponding actions of the first-base runner in those states. A mathematical treatment of policies is presented in the Supplementary Material. We denote the set of policies for the runner and pitcher by $\Pi_\text{R}$ and $\Pi_\text{P}$, respectively.

    The concept of policies allows us to mathematically characterize the runs scored in an inning. Suppose the first-base runner and the pitcher employ policies $\pi_\text{R} \in \Pi_\text{R}$ and $\pi_\text{P} \in \Pi_\text{P}$, respectively. Assuming the play starts from a state $s \in {\cal S}$, the expected total runs scored (by the batting team) in an inning then equals 
    \begin{align}
       \label{eqn:sg_pol_val}
       V^{\pi_\text{R},\pi_\text{P}}(s) &\defeq \mathbb{E}^{\pi_\text{R},\pi_\text{P}}\Bigg[\sum_{t=0}^\infty r(S_{t+1}|S_t) \Big| S_0 = s\Bigg], 
    \end{align}
    where $\mathbb{E}^{\pi_\text{R},\pi_\text{P}}$ denotes the expectation with respect to the distribution induced by $\pi_\text{R}$ and $\pi_{\text{P}}$. Observe that, unlike the typical sequential games studied in the literature, we do not discount the rewards. This is because our optimization objective is rest-of-inning expected runs, which by definition does not discount runs scored on later plays. Though the undiscounted case is more realistic in our setup, it might lead to technical issues since $V^{\pi_\text{R},\pi_\text{P}}(s)$ in \eqref{eqn:sg_pol_val} can potentially be infinite. However, all innings eventually end. As a result, the infinite series appearing inside the expectation in the right-hand side of \eqref{eqn:sg_pol_val} can be approximated with arbitrary precision by truncating the sum at some finite $T \in \mathbb{N}$, and hence we can expect $V^{\pi_\text{R},\pi_\text{P}}(s)$ to take finite values. This reasoning is rigorously justified in the Supplementary Material. Therefore, in the rest of the paper, we take $V^{\pi_\text{R},\pi_\text{P}}(s)$ to be a finite number for all $\pi_\text{R} \in \Pi_\text{R}, \pi_\text{P} \in \Pi_\text{P}$, and $s \in {\cal S}$. 
   
    The goal of the batting team is to select a policy that will maximize the expected total runs scored by them in the inning. At the same time, the fielding team minimizes the expected total runs scored by the batting team by employing a suitable policy. Taking the rational behavior of both the players into account, the maximum runs scored by the batting team when a play begins from state $s \in {\cal S}$ is given as 
    \begin{equation}
    \label{eqn:sg_val}
      V^*(s) \defeq \max_{\pi_\text{R} \in \Pi_\text{R}}\min_{\pi_\text{P} \in \Pi_\text{P}}V^{\pi_\text{R},\pi_\text{P}}(s) \ \forall s \in S.
    \end{equation}
    We call $V^*$ the \textit{optimal value function}. In our context, $V^*$ is the maximum expected total runs scored in an inning by the batting team, which, by the zero-sum property of the game, equals the minimum expected total runs incurred by the fielding team. A policy pair $(\pi_\text{R}^*,\pi_\text{P}^*) \in \Pi_\text{R} \times \Pi_\text{P}$ such that $V^{\pi_\text{R}^*,\pi_\text{P}^*}(s) = V^*(s) \ \forall s \in {\cal S}$ is called an equilibrium policy. That is, $\pi_\text{R}^*$ is a policy for the batting team that maximizes the expected total runs scored by them, and $\pi_\text{P}^*$ is a policy for the fielding team that minimizes the expected total runs incurred by them. In other words, $\pi_\text{R}^*$ is an optimal policy for the batting team and $\pi_\text{P}^*$ is an optimal policy for the fielding team. 

    A characterization of optimality conditions for \eqref{eqn:sg_val} including a dynamic programming algorithm is presented in the Supplementary Material. In our computational experiments in Section \ref{sec:results}, we use value iteration to solve \eqref{eqn:sg_val}. The value iteration algorithm is as follows. Starting from an arbitrary $V^0 \in \mathbb{R}^{\abs{\cal S}}$ with $V^0(\Delta) = 0$, iteratively obtain $V^k$ via the following updates
    \begin{align}
    \label{eqn:vi_algorithm}
      V^k(s) &\defeq \max_{a_\text{R} \in A_\text{R}(s)} \min_{a_\text{p} \in A_\text{P}(s,a_\text{R})}\sum_{s' \in \mathcal{S}}[r(s'|s) + V^{k-1}(s')] p(s'|s,a_\text{R},a_\text{P}) \quad \forall s \in \mathcal{S}.
    \end{align}
    We establish in the Supplementary Material that the iterates $V^k$ generated in \eqref{eqn:vi_algorithm} converge to the optimal value function $V^*$ at a linear rate.
    
    We note that fixing a pitcher policy reduces our two-player sequential game to a one-player MDP whose objective is to optimize the lead distance for the first-base runner. Optimality criterion similar to the two-player sequential game also exists for MDP with undiscounted rewards \parencite[Section 7.2]{bertsekas_dynamic_2012}. Numerical experiments using this decision  model for a specific pitcher policy are presented in Section \ref{sec:results-one-player}. These results may provide more immediately actionable insights for lead distance selection, as we discuss in Section \ref{sec:discussion}.

  \section{Hierarchical transition probability model}
  \label{sec:transition-probability-model}

    To infer the two-player sequential game model using empirical data, it is critical to estimate the transition probabilities $p(\cdot \mid s,a_\text{R},a_\text{P})$. In other words, we must understand how the runner's chosen lead distance and the pitcher's chosen action (pickoff or pitch) influence the stochastic transitions between game states. For states with trivial action space ($a_\text{R} = a_\text{P} = \delta$), we simply use the empirical transition probabilities between states observed in the data. For more interesting states, we decompose the state transition into (1) a runner outcome and (2) a conditional state transition, given the runner outcome. While the runner outcome depends on runner and pitcher actions, we assume that the state transition is conditionally independent of the actions, given the runner outcome. In other words, the actions only influence the state transition through the runner outcome.

    Throughout this section, we use $s_i = (b_i, c_i, d_i, o_i)$ to denote the state of the game at the beginning of each observed play $i \in \{1, ..., n\}$ (which can be a pitch or a pickoff attempt). We use $\mathcal{H}^R$, $\mathcal{H}^P$, and $\mathcal{H}^C$ to denote the sets of runners, pitchers and catchers, respectively. On play $i$:
    \begin{itemize}
      \item $h_i^R \in \mathcal{H}^R$ is the runner;
      \item $h_i^P \in \mathcal{H}^P$ is the pitcher;
      \item $h_i^C \in \mathcal{H}^C$ is the catcher;
      \item $z_i^R \in \mathbb{R^+}$ is the sprint speed of runner $h_i^R$; and
      \item $z_i^C \in \mathbb{R^+}$ is the arm strength of catcher $h_i^C$.
    \end{itemize}
    
    When the play begins, we use $\ell_i \in \mathbb{R}^+$ to denote the lead distance (in feet) taken by the first-base runner; and we use $p_i \in \{0, 1\}$ to denote whether the pitcher attempts (1) or does not attempt (0) a pickoff. We use $r_i \in \{\mbox{PO}^+,\, \mbox{PO}^-,\, \mbox{SB}^+,\, \mbox{SB}^-,\, \mbox{NG},\, \mbox{GI}\} \equiv \mathcal{R}$ to denote the runner outcome, where
    \begin{itemize}
      \item $\mbox{PO}^+$ represents a successful pickoff attempt (runner is out);
      \item $\mbox{PO}^-$ represents an unsuccessful pickoff attempt (runner is safe);
      \item $\mbox{SB}^+$ represents a successful stolen base attempt (runner is safe);
      \item $\mbox{SB}^-$ represents an unsuccessful stolen base attempt (runner is out); and
      \item $\mbox{NG}$ represents a pitch with the runner not going; and
      \item $\mbox{GI}$ represents a pitch with the runner going but interrupted by a batter event.
    \end{itemize}
    For example, a runner going on the pitch might be interrupted by a batted ball (fair or foul) or a walk causing the runner to freely advance to the next base. In these cases, the official scorebook does not acknowledge a stolen base attempt by the runner. We distinguish between NG and GI because the conditional distributions of state transitions differ between these two events. We do not define $\ell_i$ and $r_i$ for plays with no runner on first base.
    
    Using $\mathcal{E}$ to denote the event $\{S_i = s, \ell_i = a_\text{R}, p_i = a_\text{P}\}$ for states $s$ with nontrivial action space, the law of total probability decomposes the transition probability function as
    \begin{align}
      \label{eqn:prob-state-transition-two-player}
        p(s'|s,\, a_\text{R},\, a_\text{P}) = \mathbb{P}(S_{i+1} = s' \mid \mathcal{E}) = \sum_{r \in \mathcal{R}} \big(\mathbb{P}(R_i = r \mid \mathcal{E}) \cdot \mathbb{P}(S_{i+1} = s' \mid \mathcal{E},\, R_i = r)\big).
    \end{align}
    In Section \ref{sec:prob-runner-outcome}, we explain the inference of the runner outcome probability $\big(\mathbb{P}(R_i = r \mid \mathcal{E})$, and in Section \ref{sec:prob-state-transition}, we explain the inference of the conditional state transition probability $\mathbb{P}(S_{i+1} = s' \mid \mathcal{E},\, R_i = r)$.

    \subsection{Runner outcome model}
    \label{sec:prob-runner-outcome}

      There are $|\mathcal{R}| = 6$ possible outcomes for the runner with respect to the run game, and they have a hierarchical structure illustrated by Figure \ref{fig:runner-outcome-tree}. We leverage this hierarchy to decompose the runner outcome probability into five components. Using $R_i$ to denote the random variable governing the probability distribution over $r_i$,
      \begin{itemize}
        \item $\phi_i = \mathbb{P}(R_i \in \{\mbox{PO}^+,\, \mbox{PO}^-\})$ is the pickoff attempt probability;
        \item $\phi_i^+ = \mathbb{P}(R_i \in \{\mbox{PO}^+\} \mid R_i \in \{\mbox{PO}^+,\, \mbox{PO}^-\})$ is the pickoff success probability;
        \item $\psi_i = \mathbb{P}(R_i \in \{\mbox{GI},\, \mbox{SB}^+,\, \mbox{SB}^-\} \mid R_i \notin \{\mbox{PO}^+,\, \mbox{PO}^-\})$ is the runner going probability;
        \item $\nu_i = \mathbb{P}(R_i \in \{\mbox{GI}\} \mid R_i \in \{\mbox{GI}, \mbox{SB}^+,\, \mbox{SB}^-\})$ is the batter interruption probability; and
        \item $\psi_i^+ = \mathbb{P}(R_i \in \{\mbox{SB}^+\} \mid R_i \in \{\mbox{SB}^+,\, \mbox{SB}^-\})$ is the stolen base success probability.
      \end{itemize}

      \begin{figure}
        \begin{tikzpicture}
          \node (R)   at ( 0.0,  0) {$\mathcal{R}$};
          \node (PO)  at ( 2, 1) {PO};
          \node (P)   at (2, -1) {P};
          \node (PO+) at ( 4.0, 1.5) {PO$^+$};
          \node (PO-) at ( 4, 0.5) {PO$^-$};
          \node (G)   at (4, -0.5) {G};
          \node (NG)  at (4, -1.5) {NG};
          \node (GI)  at (6, 0) {GI};
          \node (SB)  at (6, -1) {SB};
          \node (SB+) at (8, -0.5) {SB$^+$};
          \node (SB-) at (8, -1.5) {SB$^-$};
          \draw[->] (R)  -- node[above] {\scriptsize $\phi_i$} (PO);
          \draw[->] (R)  -- (P);
          \draw[->] (PO) -- node[above] {\scriptsize $\phi_i^+$} (PO+);
          \draw[->] (PO) -- (PO-);
          \draw[->] (P)  -- node[above] {\scriptsize $\psi_i$} (G);
          \draw[->] (P)  -- (NG);
          \draw[->] (G)  -- node[above] {\scriptsize $\nu_i$} (GI);
          \draw[->] (G)  -- (SB);
          \draw[->] (SB) -- node[above] {\scriptsize $\psi_i^+$} (SB+);
          \draw[->] (SB) -- (SB-);
        \end{tikzpicture}
        \caption{
          \it A diagram of the runner outcome hierarchy. Each upward-moving edge is labeled with the conditional probability of the top child node given the parent node. The conditional probability of the bottom child is one minus the conditional probability of the top child. The sets PO (pickoff attempt), P (pitch), G (runner going) and SB (stolen base attempt) are subsets of the runner outcome set $\mathcal R$.
        }
        \label{fig:runner-outcome-tree}
      \end{figure}

      Each of these probabilities depends on the play's contextual covariates and player skills, for which the subscript $i$ accounts. If the pitcher has agency, then the pickoff attempt probability depends only on pitcher action, $\phi_i = p_i$. However, it is useful for some applications to remove the pitcher's agency, as described in Section \ref{sec:prob-pickoff-attempt}, so we keep the notation general. The runner outcome probability comprises the products of the component probabilities:
      \begin{align}
        \label{eqn:prob-runner-outcome-two-player}
        \begin{split}
          \mathbb{P}(R_i = r \mid \mathcal{E}) =
          \begin{cases}
            \phi_i \cdot \phi_i^+                                             & \mbox{if } r = \mbox{PO}^+\\
            \phi_i \cdot (1 - \phi_i^+)                                       & \mbox{if } r = \mbox{PO}^-\\
            (1 - \phi_i) \cdot \psi_i \cdot \nu_i                             & \mbox{if } r = \mbox{GI}\\
            (1 - \phi_i) \cdot \psi_i \cdot (1 - \nu_i) \cdot \psi_i^+        & \mbox{if } r = \mbox{SB}^+\\
            (1 - \phi_i) \cdot \psi_i \cdot (1 - \nu_i) \cdot (1 - \psi_i^+)  & \mbox{if } r = \mbox{SB}^-\\
            (1 - \phi_i) \cdot (1 - \psi_i)                                   & \mbox{if } r = \mbox{NG}\\
          \end{cases}.
        \end{split}
      \end{align}

      We estimate the component binomial outcome probabilities using GLMMs \parencite{stroup_generalized_2012}. Mixed-effects models allow us to incorporate player random effects along with fixed effects for game state variables and lead distance. Player random effects are particularly appropriate for this setting because they yield player skill estimates which go beyond measurable characteristics like sprint speed and are regularized towards a common mean. Sections \ref{sec:prob-pickoff-attempt} through \ref{sec:prob-stolen-base} detail the model specifications.

      \subsubsection{Pickoff attempt probability}
      \label{sec:prob-pickoff-attempt}

        In the two-player sequential game model, there is no need to estimate pickoff attempt probability because it is a deterministic decision made by the pitcher. However, as we discuss in Section \ref{sec:discussion}, modeling the pickoff decision as stochastic rather than antagonistic results in different, and perhaps more immediately actionable, insights for the runner. From the runner's perspective, this approach models the pitcher's action as a stochastic response to the runner's lead distance. In this section we present a probability model for pickoff attempts, which forms the basis for the results in Section \ref{sec:results-one-player}.

        We model $\phi_i$ using a mixed-effects logistic regression with fixed effects for balls, strikes, outs, disengagements (all as categorical), and lead distance (numeric); and with a random effect for the pitcher.
        \begin{align}
          \label{eqn:prob-pickoff-attempt}
          \begin{split}
            \log\left(\frac{\phi_i}{1 - \phi_i}\right) &= \alpha + \beta^B_{c_i^B} + \beta^S_{c_i^S} +
              \beta^O_{o_i} + \beta^D_{d_i} + \beta^L \ell_i + \gamma^P_{h_i^P}\\
            \gamma^P_{h} &\sim \mathcal{N}(0, \sigma^2_P) \hspace{4mm} \mbox{\it i.i.d.} \hspace{4mm} \forall h \in \mathcal{H}^P.
          \end{split}
        \end{align}
        This model has eleven fixed, unknown parameters: the intercept $\alpha$; three $\beta^B$'s, two $\beta^S$'s, two $\beta^O$'s, two $\beta^D$'s, and $\beta^L$; and the variance parameter $\sigma^2_P$. The number of random effects is $|\mathcal{H}^P|$, and we interpret them as latent individual tendencies to attempt pickoffs.
        
      \subsubsection{Pickoff success probability}
      \label{sec:prob-pickoff-success}

        We model $\phi_i^+$ using a mixed-effects logistic regression with a fixed effect for lead distance and a random effect for pitcher.
        \begin{align}
          \label{eqn:prob-pickoff-success}
          \begin{split}
            \log\left(\frac{\phi_i^+}{1 - \phi_i^+}\right) &= \alpha + \beta^L \ell_i + \gamma^P_{h_i^P},\\
            \gamma^P_{h} &\sim \mathcal{N}(0, \sigma^2_P) \hspace{4mm} \mbox{\it i.i.d.} \hspace{4mm} \forall h \in \mathcal{H}^P.
          \end{split}
        \end{align}
        This model has three fixed, unknown parameters: the intercept $\alpha$, the slope $\beta^L$, and the variance parameter $\sigma^2_P$. The number of random effects is $|\mathcal{H}^P|$, and we interpret them as latent individual pickoff skill.

      \subsubsection{Runner going probability}
      \label{sec:prob-runner-going}

        Modeling the runner's decision to attempt a stolen base presents empirical challenges. While classical game theory might frame this as a deterministic decision, this model aligns poorly with real-life decisions made by baserunners. We posit that the decision is a response to latent environment variables (e.g., the pitcher's body language or timing between pitches) that we do not observe in the data available. Therefore, we treat the decision as stochastic given player skills and contextual factors. Further, we exclude lead distance from this model to avoid a paradoxical situation in which the optimal action for the runner would be to shorten their lead to reduce the probability of an ill-advised stolen base attempt.

        We model $\psi_i$ using a mixed-effects logistic regression with fixed effects for balls, strikes, outs, disengagements, runner sprint speed (numeric), and catcher arm strength (numeric); and with random effects for runner, pitcher, and catcher.
        \begin{align}
          \label{eqn:prob-runner-going}
          \begin{split}
            \log\left(\frac{\psi_i}{1 - \psi_i}\right) &= \alpha + \beta^B_{c_i^B} + \beta^S_{c_i^S} +
              \beta^O_{o_i} + \beta^D_{d_i} + (\beta^R z_i^R + \gamma^R_{h_i^R}) +
              \gamma^P_{h_i^P} + (\beta^C z_i^C + \gamma^C_{h_i^C}),\\
            \gamma^R_{h} &\sim \mathcal{N}(0, \sigma^2_R) \hspace{4mm} \mbox{\it i.i.d.} \hspace{4mm} \forall h \in \mathcal{H}^R,\\
            \gamma^P_{h} &\sim \mathcal{N}(0, \sigma^2_P) \hspace{4mm} \mbox{\it i.i.d.} \hspace{4mm} \forall h \in \mathcal{H}^P,\\
            \gamma^C_{h} &\sim \mathcal{N}(0, \sigma^2_C) \hspace{4mm} \mbox{\it i.i.d.} \hspace{4mm} \forall h \in \mathcal{H}^C.
          \end{split}
        \end{align}
 
        We interpret $(\beta^R z_i^R + \gamma^R_{h_i^R})$ as the effect of runner $h_i^R$ on the runner going probability, combining the effect of the runner's measured sprint speed with a latent effect representing the individual's tendency to attempt a steal. By including a fixed effect for sprint speed, we are effectively regularizing the runner effects toward priors based on their sprint speeds. The interpretation of $(\beta^C z_i^C + \gamma^C_{h_i^C})$ for catchers is similar: Each catcher has an estimated effect regularized toward a prior based on their arm strength.
      
      \subsubsection{Batter interruption probability}      \label{sec:prob-going-interrupt}

        We model $\nu_i$ using a mixed-effects logistic regression with fixed effects for balls, strikes, outs, disengagements, runner sprint speed, and catcher arm strength; and with random effects for runner, pitcher, and catcher.
        \begin{align}
          \label{eqn:prob-going-interrupt}
          \begin{split}
            \log\left(\frac{\nu_i}{1 - \nu_i}\right) &= \alpha + \beta^B_{c_i^B} + \beta^S_{c_i^S} +
              \beta^O_{o_i} + \beta^D_{d_i} + (\beta^R z_i^R + \gamma^R_{h_i^R}) +
              \gamma^P_{h_i^P} + (\beta^C z_i^C + \gamma^C_{h_i^C}),\\
            \gamma^R_{h} &\sim \mathcal{N}(0, \sigma^2_R) \hspace{4mm} \mbox{\it i.i.d.} \hspace{4mm} \forall h \in \mathcal{H}^R,\\
            \gamma^P_{h} &\sim \mathcal{N}(0, \sigma^2_P) \hspace{4mm} \mbox{\it i.i.d.} \hspace{4mm} \forall h \in \mathcal{H}^P,\\
            \gamma^C_{h} &\sim \mathcal{N}(0, \sigma^2_C) \hspace{4mm} \mbox{\it i.i.d.} \hspace{4mm} \forall h \in \mathcal{H}^C.
          \end{split}
        \end{align}
      
      \subsubsection{Stolen base success probability}
      \label{sec:prob-stolen-base}

        We model $\psi_i^+$ using a mixed-effects logistic regression with fixed effects for lead distance, runner sprint speed, and catcher arm strength; and with random effects for runner, pitcher, and catcher.
        \begin{align}
          \label{eqn:prob-stolen-base}
          \begin{split}
            \log\left(\frac{\psi_i^+}{1 - \psi_i^+}\right) &= \alpha + \beta^L \ell_i + (\beta^R z_i^R + \gamma^R_{h_i^R}) +
              \gamma^P_{h_i^P} + (\beta^C z_i^C + \gamma^C_{h_i^C}),\\
            \gamma^R_{h} &\sim \mathcal{N}(0, \sigma^2_R) \hspace{4mm} \mbox{\it i.i.d.} \hspace{4mm} \forall h \in \mathcal{H}^R,\\
            \gamma^P_{h} &\sim \mathcal{N}(0, \sigma^2_P) \hspace{4mm} \mbox{\it i.i.d.} \hspace{4mm} \forall h \in \mathcal{H}^P,\\
            \gamma^C_{h} &\sim \mathcal{N}(0, \sigma^2_C) \hspace{4mm} \mbox{\it i.i.d.} \hspace{4mm} \forall h \in \mathcal{H}^C.
          \end{split}
        \end{align}
        As in (\ref{eqn:prob-runner-going}), we interpret $(\beta^R z_i^R + \gamma^R_{h_i^R})$ as the runner effect regularized toward a prior based on sprint speed and $(\beta^C z_i^C + \gamma^C_{h_i^C})$ as the catcher effect regularized toward a prior based on arm strength. Through these terms we account for both measurable skills and latent skills that influence the probability of a successful stolen base.

    \subsection{Conditional state transition model}
    \label{sec:prob-state-transition}

      Using the assumption that state transitions depend on the runner and pitcher actions only through the runner outcome, we simplify the conditional state transition probability as
      \begin{align}
        \mathbb{P}(S_{i+1} = s' \mid \mathcal{E},\, R_i = r) =
          \mathbb{P}(S_{i+1} = s' \mid S_i = s,\, R_i = r).
      \end{align}
      To estimate this probability, we use conditional empirical frequencies of the ending state $s'$ given the starting state $s$ and the runner outcome $r$. The disengagements $d$ present a problem because few plays begin with $d = 1$ and even fewer begin with $d = 2$. However, there is little reason to expect that $d$ influences state transitions beyond its influence on the runner outcome. Therefore, we pool across $d$ when calculating empirical frequencies.
      
      In detail, we define reduced states $\tilde s \equiv (b, c, \cdot, o)$ and $\tilde s' \equiv (b', c', \cdot, o')$. The pooled empirical transition probability $Q$ conditions on $e = d' - d$ instead of conditioning on $d$ and $d'$. Using $E_{i}$ to denote the random variable governing the probability distribution over $d_{i + 1} - d_i$,
      \begin{align}
        \label{eqn:pooled-probabilities}
          Q(\tilde s, r, \tilde s', e)
            = \hat{\mathbb{P}}\left(\tilde S_{i+1} =   \tilde s',\, E_i = e \mid \tilde S_i = \tilde s,\, R_i = r\right)
            = \frac
              {\sum_{i=1}^n 1_i^{\tilde s, r} \cdot 1_i^{\tilde s', e}}
              {\sum_{i=1}^n 1_i^{\tilde s, r}},
      \end{align}
      where $1_i^{\tilde s, r} = \mathbb{I}\left\{\tilde s_i = \tilde s,\, r_i = r\right\}$ and $1_i^{\tilde s', e}= \mathbb{I}\left\{\tilde s_{i+1} = \tilde s',\, d_{i+1} - d_i = e\right\}$. Finally, we estimate the full-state conditional transition probabilities as
      \begin{align*}
        &\hat{\mathbb{P}}(S_{i+1} = (b', c', d', o') \mid S_i = (b, c, d, o), R_i = r) = Q\big((b, c, \cdot, o),\, r,\, (b', c', \cdot, o'),\, d' - d\big).
      \end{align*}
  
    \subsection{Validation and uncertainty quantification}

      We perform 20-fold cross-validation to estimate the out-of-sample prediction error, using negative log loss (deviance). We report the performance as the percentage of deviance explained relative to the null model which always predicts the global mean probability of a positive outcome. A low percentage of deviance explained does not necessarily indicate poor model fit. It would be unrealistic to expect the models to predict individual binary outcomes with high certainty. The purpose of these models is to reflect how player identities, lead distance, and context predict subtle shifts in outcome probabilities which influence the optimal lead distance policy.

      The end product of our analysis is the optimal policy estimated via value iteration as described in Section \ref{sec:sg-model}. The optimal policy is a function of the transition probabilities estimated in this section. We propagate the uncertainty in the transition probability estimation through to the end result using a bootstrapping procedure. For each of $B = 100$ bootstrap replicates, we refit the GLMMs from Section \ref{sec:prob-runner-outcome} via parametric bootstrapping \parencite{bates_fitting_2015} and re-estimate the conditional transition probabilities from Section \ref{sec:prob-state-transition} by resampling the state transitions in the original dataset, stratified by year. We solve the optimal policy for each bootstrap replicate and report the uncertainty as the standard deviation of the end result across the replicates.

  \section{Results}
  \label{sec:results}

    In Section \ref{sec:results-multilevel-model} we present the results of estimating the GLMMs described in Section \ref{sec:prob-runner-outcome}, and in Section \ref{sec:results-two-player} we present the results of solving the two-player sequential game. Finally, in Section \ref{sec:results-one-player}, we consider a reduction of the two-player game by removing the pitcher's agency and replacing it with a stochastic model. With this substitution, the two-player game becomes a one-player MDP, for which we find the optimal policy.
    
    \subsection{Runner outcome probability model}
    \label{sec:results-multilevel-model}

\begin{figure}
        \centering
        \includegraphics[width = 0.4\linewidth]{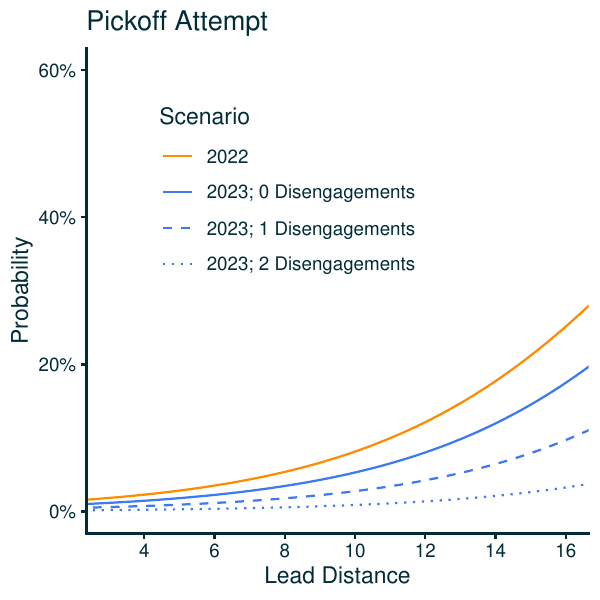}
        \includegraphics[width = 0.4\linewidth]{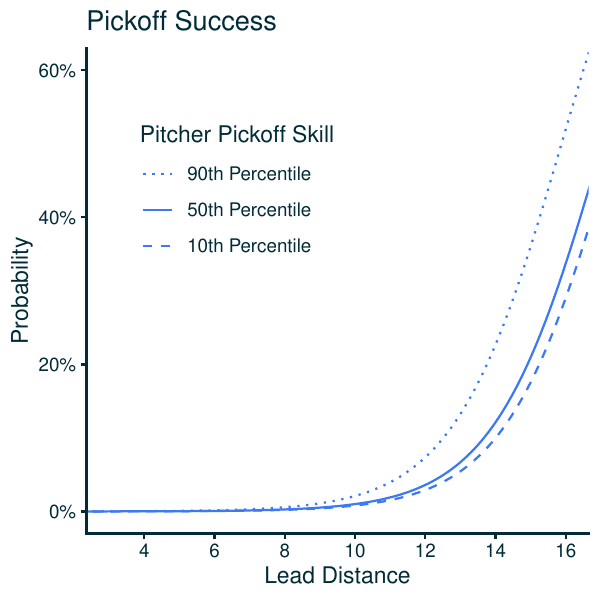}
        \caption{
          \it Pickoff attempt probability (left) and success probability (right) modeled as functions of lead distance. The left figure corresponds to the model detailed from Section \ref{sec:prob-pickoff-attempt}, assuming zero balls, zero strikes, zero outs, and a median pitcher random effect. The right figure corresponds to the model from Section \ref{sec:prob-pickoff-success}, which includes only an intercept, a fixed effect for lead distance and a random effect for pitcher identity.
        }
        \label{fig:prob-pickoff}
      \end{figure}

\begin{figure}
        \centering
        \includegraphics[width = 0.8\linewidth]{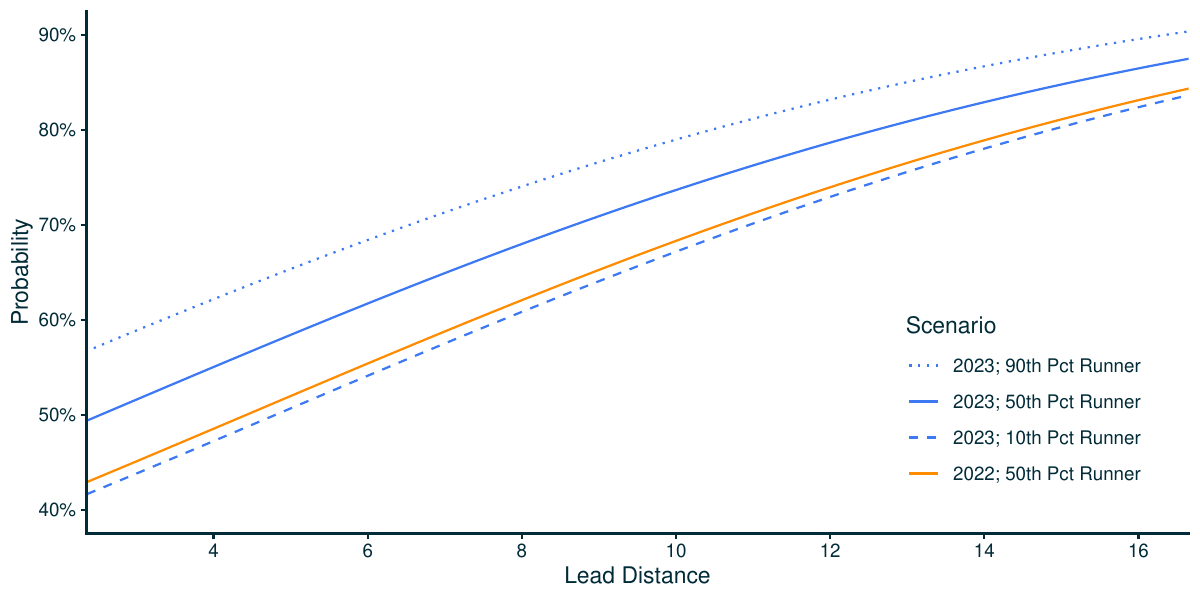}
        \caption{
          \it Stolen base success probability modeled as function of lead distance, split by season and runner effect (combining the fixed effect of runner sprint speed and the random effect of runner identity). This figure shows the probability estimated by the model detailed in Section \ref{sec:prob-stolen-base}, assuming zero balls, zero strikes, zero outs, and a median value for the combined effect of pitcher identity (random effect), catcher arm strength (fixed effect), and catcher identity (random effect).
        }
        \label{fig:prob-stolen-base}
      \end{figure}

      We implemented the GLMMs from Section \ref{sec:prob-runner-outcome} in R using the package glmmTMB \parencite{brooks_glmmtmb_2017}. Table \ref{tab:model-summary} summarizes the results of estimating the runner outcome probability models detailed in Section \ref{sec:prob-runner-outcome}. Figures \ref{fig:prob-pickoff} and \ref{fig:prob-stolen-base} illustrate the relationship between lead distance and runner outcome probabilities for the models which include lead distance. For each of these three models, the direction of the relationship between lead distance and outcome probability is as expected. In 20-fold cross-validation, the component GLMMs explained 5.4\% of deviance in pickoff attempt probability; 9.9\% of deviance in pickoff success probability; 16.7\% of deviance in runner going probability; 13.0\% of deviance in batter interruption probability; and 2.5\% of deviance in stolen base success probability.

      \begin{table}
        \resizebox{\linewidth}{!}{
          \begin{tabular}{l|rrrrr}
&\multicolumn{5}{c}{Estimated Effect on Log-Odds}\\
Variable & Pickoff Attempt & Pickoff Success & Runner Going & Batter Interruption & Stolen Base \\
  \hline
{\it Intercept} & $-2.42$ {\color{gray} $\pm0.03$} & $-5.27$ {\color{gray} $\pm0.15$} & $-4.87$ {\color{gray} $\pm0.06$} & $-0.87$ {\color{gray} $\pm0.08$} & $0.68$ {\color{gray} $\pm0.08$} \\ 
   \hline
{\it Fixed Effects} &  &  &  &  &  \\ 
  Year (2023 vs. 2022) & $-0.46$ {\color{gray} $\pm0.02$} & $0.75$ {\color{gray} $\pm0.13$} & $-0.11$ {\color{gray} $\pm0.03$} & $-0.09$ {\color{gray} $\pm0.06$} & $0.26$ {\color{gray} $\pm0.08$} \\ 
  Outs (1 vs. 0) & $0.09$ {\color{gray} $\pm0.02$} &  & $0.27$ {\color{gray} $\pm0.03$} & $-0.10$ {\color{gray} $\pm0.06$} &  \\ 
  Outs (2 vs. 0) & $-0.07$ {\color{gray} $\pm0.02$} &  & $0.46$ {\color{gray} $\pm0.03$} & $0.14$ {\color{gray} $\pm0.06$} &  \\ 
  Balls (1 vs. 0) & $-0.04$ {\color{gray} $\pm0.02$} &  & $0.17$ {\color{gray} $\pm0.03$} & $0.11$ {\color{gray} $\pm0.06$} &  \\ 
  Balls (2 vs. 0) & $-0.35$ {\color{gray} $\pm0.03$} &  & $-0.02$ {\color{gray} $\pm0.04$} & $0.23$ {\color{gray} $\pm0.07$} &  \\ 
  Balls (3 vs. 0) & $-0.41$ {\color{gray} $\pm0.04$} &  & $1.70$ {\color{gray} $\pm0.04$} & $1.57$ {\color{gray} $\pm0.08$} &  \\ 
  Strikes (1 vs. 0) & $0.05$ {\color{gray} $\pm0.02$} &  & $0.19$ {\color{gray} $\pm0.03$} & $0.54$ {\color{gray} $\pm0.06$} &  \\ 
  Strikes (2 vs. 0) & $0.03$ {\color{gray} $\pm0.02$} &  & $0.86$ {\color{gray} $\pm0.03$} & $1.13$ {\color{gray} $\pm0.07$} &  \\ 
  Disengagements (1 vs. 0) & $-0.68$ {\color{gray} $\pm0.04$} &  & $0.84$ {\color{gray} $\pm0.04$} & $-0.15$ {\color{gray} $\pm0.07$} &  \\ 
  Disengagements (2 vs. 0) & $-1.84$ {\color{gray} $\pm0.23$} &  & $1.54$ {\color{gray} $\pm0.08$} & $-0.06$ {\color{gray} $\pm0.13$} &  \\ 
  Lead Distance (ft) & $0.22$ {\color{gray} $\pm0.00$} & $0.65$ {\color{gray} $\pm0.05$} &  &  & $0.14$ {\color{gray} $\pm0.03$} \\ 
  Runner Sprint Speed (ft/s) &  &  & $0.60$ {\color{gray} $\pm0.02$} & $-0.15$ {\color{gray} $\pm0.02$} & $0.18$ {\color{gray} $\pm0.04$} \\ 
  Catcher Arm Strength (mi/h) &  &  & $-0.02$ {\color{gray} $\pm0.01$} & $0.01$ {\color{gray} $\pm0.01$} & $-0.06$ {\color{gray} $\pm0.01$} \\ 
   \hline
{\it Random Effects} &  &  &  &  &  \\ 
  Pitcher &  {\color{gray} $\pm0.60$} &  {\color{gray} $\pm0.93$} &  {\color{gray} $\pm0.58$} &  {\color{gray} $\pm0.17$} &  {\color{gray} $\pm0.30$} \\ 
  Runner &  &  &  {\color{gray} $\pm0.67$} &  {\color{gray} $\pm0.22$} &  {\color{gray} $\pm0.29$} \\ 
  Catcher &  &  &  {\color{gray} $\pm0.17$} &  {\color{gray} $\pm0.11$} &  {\color{gray} $\pm0.23$} \\ 
\end{tabular}

        }
        \caption{
          \it Summary of estimated runner outcome probability models. We estimated four different generalized (binomial) linear mixed-effects models for pickoff (PO) attempt and success and for stolen base (SB) attempt and success. For fixed effects we report the estimated effect with one standard error, and for random effects we report the estimated standard deviation term. For 2022, we coded disengagements as always zero to reflect no pickoff limit.
        }
        \label{tab:model-summary}
       \end{table}
      
      Figure \ref{fig:prob-pickoff} illustrates how pickoff attempt probability and success probability vary with distance. Pickoff attempt rates were lower in 2023 than in 2022, especially as the number of prior disengagements increased. The modeled probability remains low even for extremely long leads. The ball-strike count and the identity of the pitcher also influence the pickoff attempt probability, but the primary drivers are lead distance and the number prior disengagements. Except for the most skilled pitchers, pickoff success probability is very low within the range covering the vast majority of lead distances. Toward the high end of this range, the success rate increases quickly as lead distance increases. Interestingly, there is a large gap between the 90$^\text{th}$-percentile pitcher and the median pitcher but a small gap between the median pitcher and the 10$^\text{th}$-percentile pitcher, demonstrating a right skew in the distribution of pitcher pickoff skill.
 
      Figure \ref{fig:prob-stolen-base} shows stolen base success probability modeled as a function of lead distance, split by season and by runner skill. Conditioned on lead distance, 10$^\text{th}$-percentile runners in 2023 had similar success probabilities to the median runner in 2022, which is likely attributable to rule changes making it easier to steal bases in 2023. Base-out state and battery (i.e. pitcher and catcher) skill also influence stolen base success probability although this visualization does not show them.
      
      A novel applied result of these runner outcome probability models stems from the combination of fixed effects for measurable player characteristics (catcher arm strength and runner sprint speed) and random effects for player identities. This feature allows us to answer questions such as: How important are these measurable characteristics for player performance outcomes? And to what extent do players over- or under-perform what is expected from their measurable characteristics? The combination of the fixed effect and the random effect gives us an estimate (which we call the combined effect) of each player's skill, regularized toward what is expected from their measurable characteristics.

\begin{figure}
        \centering
        \includegraphics[width = 0.4\linewidth]{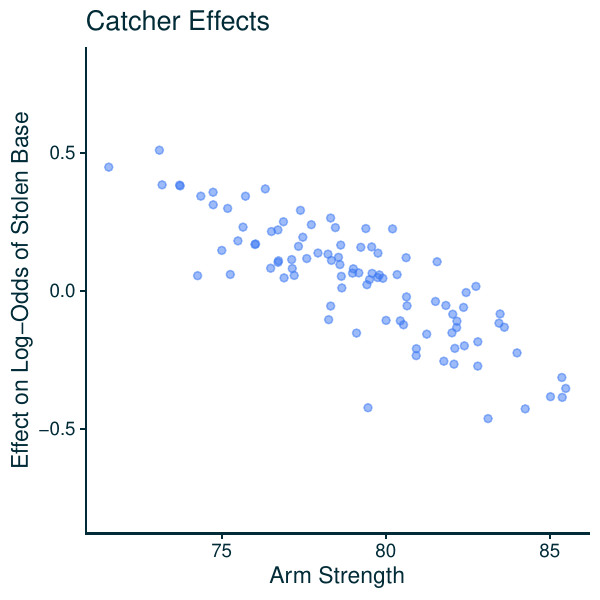}
        \includegraphics[width = 0.4\linewidth]{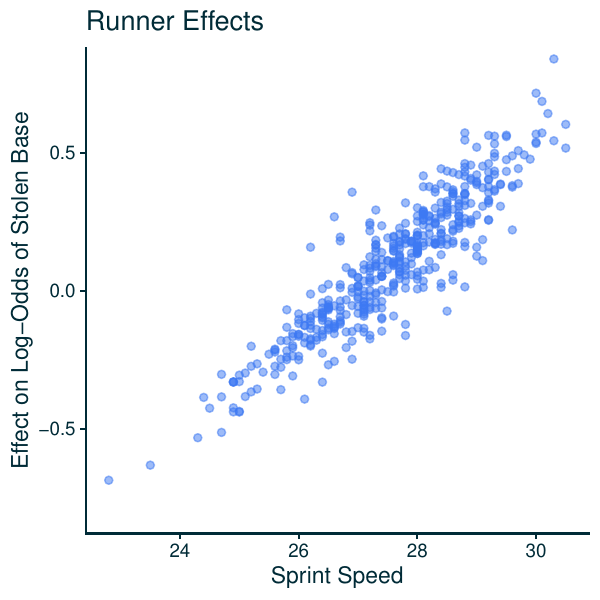}
        \caption{
          \it
          (Left) The relationship between arm strength and catcher combined effect on stolen base success probability (relative to average, on the log-odds scale). The combined effect is the sum of the fixed effect due to arm strength and the catcher random effect.
          (Right) The relationship between sprint speed and runner combined effect on stolen base success probability (relative to average, on the log-odds scale). The combined effect is the sum of the fixed effect due to sprint speed and the runner random effect.
        }
        \label{fig:random-effect}
      \end{figure}

      For the stolen base success probability model, Figure \ref{fig:random-effect} shows the relationship between the combined player effect and the corresponding measurable characteristic, for both catchers and runners. Arm strength explains 70\% of the variance in the catcher combined effect and that sprint speed explains 84\% of the variance in the runner combined effect. In other words, sprint speed tells us more about runner skill than arm strength tells us about catcher skill, with regard to stolen base outcomes. Beyond arm strength, the catcher can reduce stolen base success by releasing the ball more quickly or by throwing it more accurately to the target base. Beyond sprint speed, the runner can increase stolen base success by starting to run earlier (i.e., getting a better jump) or by accelerating more quickly to their top speed.
      
    \subsection{Two-player sequential game}
    \label{sec:results-two-player}
    
      Using the transition probabilities estimated on the empirical data, we implemented the value iteration algorithm from Section \ref{sec:sg-model}, discretizing the runner's lead distance action space to the nearest tenth of a foot. Table \ref{tab:policy-by-count-zsg} reports how the optimal runner policy varies by count and disengagements, and Table \ref{tab:policy-by-outs-zsg} reports how the optimal runner policy varies by outs and disengagements. This runner policy maximizes runs scored while taking into account the rational behavior of both players, as defined by (\ref{eqn:sg_val}). For this runner policy, the pitcher is indifferent to attempting a pickoff or not.
      
      \begin{table}
        \begin{tabular}{c|ccc}
& \multicolumn{3}{c}{Optimal Lead by Prior Disengagements}\\
Count & 0 & 1 & 2 \\
  \hline
0-0 & 9.6 {\color{gray} $\pm$ 0.4} & 11.1 {\color{gray} $\pm$ 0.3} & 15.3 {\color{gray} $\pm$ 0.4} \\ 
  1-0 & 9.7 {\color{gray} $\pm$ 0.5} & 11.0 {\color{gray} $\pm$ 0.4} & 15.2 {\color{gray} $\pm$ 0.4} \\ 
  2-0 & 9.8 {\color{gray} $\pm$ 0.7} & 11.0 {\color{gray} $\pm$ 0.4} & 14.8 {\color{gray} $\pm$ 0.4} \\ 
  3-0 & 9.6 {\color{gray} $\pm$ 1.5} & 11.1 {\color{gray} $\pm$ 0.6} & 13.8 {\color{gray} $\pm$ 0.4} \\ 
   \hline
0-1 & 9.4 {\color{gray} $\pm$ 0.4} & 10.9 {\color{gray} $\pm$ 0.3} & 15.3 {\color{gray} $\pm$ 0.4} \\ 
  1-1 & 9.9 {\color{gray} $\pm$ 0.6} & 11.1 {\color{gray} $\pm$ 0.4} & 15.2 {\color{gray} $\pm$ 0.4} \\ 
  2-1 & 9.2 {\color{gray} $\pm$ 0.7} & 10.7 {\color{gray} $\pm$ 0.4} & 14.9 {\color{gray} $\pm$ 0.4} \\ 
  3-1 & 9.9 {\color{gray} $\pm$ 1.1} & 11.1 {\color{gray} $\pm$ 0.6} & 14.4 {\color{gray} $\pm$ 0.4} \\ 
   \hline
0-2 & 9.9 {\color{gray} $\pm$ 0.5} & 11.1 {\color{gray} $\pm$ 0.4} & 15.4 {\color{gray} $\pm$ 0.4} \\ 
  1-2 & 9.5 {\color{gray} $\pm$ 1.6} & 11.0 {\color{gray} $\pm$ 0.5} & 15.4 {\color{gray} $\pm$ 0.4} \\ 
  2-2 & 10.4 {\color{gray} $\pm$ 0.6} & 11.3 {\color{gray} $\pm$ 0.4} & 15.0 {\color{gray} $\pm$ 0.4} \\ 
  3-2 & 10.5 {\color{gray} $\pm$ 0.7} & 11.5 {\color{gray} $\pm$ 0.5} & 14.4 {\color{gray} $\pm$ 0.4} \\ 
\end{tabular}

        \caption{
          \it The optimal lead distance (in feet) by count and prior disengagements under the two-player sequential game, assuming an average runner on first base facing an average battery (pitcher and catcher) with no outs and third base empty. Estimates include standard errors estimated via bootstrapping.
        }
        \label{tab:policy-by-count-zsg}
      \end{table}

      The most striking observation is that the optimal lead distance is significantly longer with two prior disengagements (13.8--15.4 feet) than in other states. These recommended lead distances are near the edge of observed lead distances in the dataset but not outside of the observed range. Overall, the optimal runner policy is more aggressive than observed behavior. Recall from Section \ref{sec:introduction} that the average lead distance after zero, one, and two prior disengagements was 9.5 feet, 10.3 feet, and 11.0 feet, respectively.

      \begin{table}
        \small
        \centering       
        \begin{tabular}{c|ccc}
& \multicolumn{3}{c}{Optimal Lead by Prior Disengagements}\\
Outs & 0 & 1 & 2 \\
  \hline
  0 & 9.6 {\color{gray} $\pm$ 0.4} & 11.1 {\color{gray} $\pm$ 0.3} & 15.3 {\color{gray} $\pm$ 0.4} \\ 
    1 & 9.7 {\color{gray} $\pm$ 0.5} & 11.3 {\color{gray} $\pm$ 0.3} & 15.2 {\color{gray} $\pm$ 0.4} \\ 
    2 & 10.1 {\color{gray} $\pm$ 0.5} & 11.4 {\color{gray} $\pm$ 0.3} & 15.3 {\color{gray} $\pm$ 0.4} \\ 
\end{tabular}

        \caption{
          \it The optimal lead distance (in feet) by outs and prior disengagements under the two-player sequential game, assuming an average runner on first base facing an average battery (pitcher and catcher) with an 0-0 count and third base empty. Estimates include standard errors estimated via bootstrapping.
        }
        \label{tab:policy-by-outs-zsg}
      \end{table}

      For two-disengagement states, the optimal lead distance generally, although not strictly, shortens as the count deepens (i.e., as the numbers of balls and strikes increase). This makes intuitive sense because the plate appearance is closer to ending, which relieves the pressure on the pitcher to avoid attempting a pickoff. The lead distance generally, although not strictly, increases as the number of outs increases. This makes intuitive sense because the runner is less likely to score without more aggressive action. In other words, the riskier lead distance is more palatable because the runner has a lower baseline probability of scoring from first base with two outs. The non-monotonicity of the optimal lead distance policy may be a consequence of estimating transition probabilities using empirical data.

      Before the first disengagement in 2023, runner behavior was reasonably calibrated relative to the equilibrium. Runners exceeded the optimal lead (for an average runner) in 46\% of opportunities. However, after the first disengagement runners exceeded the recommended lead in 25\% of the time, and after the second disengagement they exceeded the recommended lead in only 1.5\% of the time. That is, the equilibrium policy would be more aggressive than observed runner behavior following unsuccessful pickoff attempts in 2023.

      Pitchers also played less aggressively in 2023 than the game-theoretic equilibrium would suggest. With zero prior disengagements, pitchers attempted pickoffs 3.7\% of the time the runner took a shorter lead than recommended and only 7.7\% of the time the runner took a longer lead than recommended. After the first disengagement, pitchers only attempted pickoffs in 4.7\% of opportunities when runners exceeded the equilibrium lead distance. Under the game-theoretic equilibrium, both runners and pitchers would play more aggressively.

    \subsection{One-player MDP}
    \label{sec:results-one-player}

      The results in Tables \ref{tab:policy-by-count-zsg} and \ref{tab:policy-by-outs-zsg} provide guidance to runners for choosing lead distances that are unexploitable from a game theory perspective. However, pitchers may be unlikely to adapt their strategy to the optimal policy in the short term. Runner actions are already observable by the defense, and current pitcher behavior seems more conservative than optimal, as noted in Section \ref{sec:results-two-player}. After all, the pitcher has the additional cognitive burden of contending with the batter. In light of this, runners may achieve better results by choosing a policy to exploit observed pitcher behavior.

      We consider replacing the pitcher's agency in the two-player sequential game with the probability model described in Section \ref{sec:prob-pickoff-attempt}. The pickoff decision becomes a stochastic response to the runner's lead distance, resulting in a one-player MDP. As with the two-player game, we implement value iteration to find the optimal runner policy given observed pitcher behavior. Table \ref{tab:policy-by-count-mdp} reports how the optimal policy varies by count and prior disengagements.

      \begin{table}
        \begin{tabular}{c|ccc}
& \multicolumn{3}{c}{Optimal Lead by Prior Disengagements}\\
Count & 0 & 1 & 2 \\
  \hline
0-0 & 9.4 {\color{gray}$\pm$ 0.3} & 11.0 {\color{gray}$\pm$ 0.3} & 13.6 {\color{gray}$\pm$ 0.3} \\ 
  1-0 & 9.5 {\color{gray}$\pm$ 0.3} & 11.2 {\color{gray}$\pm$ 0.3} & 13.6 {\color{gray}$\pm$ 0.3} \\ 
  2-0 & 9.5 {\color{gray}$\pm$ 0.3} & 11.1 {\color{gray}$\pm$ 0.3} & 13.5 {\color{gray}$\pm$ 0.4} \\ 
  3-0 & 10.5 {\color{gray}$\pm$ 0.4} & 12.1 {\color{gray}$\pm$ 0.4} & 14.1 {\color{gray}$\pm$ 0.4} \\ 
   \hline
0-1 & 9.3 {\color{gray}$\pm$ 0.3} & 11.0 {\color{gray}$\pm$ 0.3} & 13.6 {\color{gray}$\pm$ 0.3} \\ 
  1-1 & 9.5 {\color{gray}$\pm$ 0.3} & 11.1 {\color{gray}$\pm$ 0.3} & 13.6 {\color{gray}$\pm$ 0.3} \\ 
  2-1 & 9.4 {\color{gray}$\pm$ 0.3} & 11.1 {\color{gray}$\pm$ 0.3} & 13.6 {\color{gray}$\pm$ 0.4} \\ 
  3-1 & 10.2 {\color{gray}$\pm$ 0.3} & 11.8 {\color{gray}$\pm$ 0.3} & 13.9 {\color{gray}$\pm$ 0.4} \\ 
   \hline
0-2 & 9.7 {\color{gray}$\pm$ 0.3} & 11.3 {\color{gray}$\pm$ 0.3} & 13.8 {\color{gray}$\pm$ 0.3} \\ 
  1-2 & 9.7 {\color{gray}$\pm$ 0.3} & 11.4 {\color{gray}$\pm$ 0.3} & 13.9 {\color{gray}$\pm$ 0.4} \\ 
  2-2 & 9.9 {\color{gray}$\pm$ 0.3} & 11.5 {\color{gray}$\pm$ 0.3} & 13.8 {\color{gray}$\pm$ 0.4} \\ 
  3-2 & 10.2 {\color{gray}$\pm$ 0.3} & 11.7 {\color{gray}$\pm$ 0.3} & 13.7 {\color{gray}$\pm$ 0.4} \\ 
\end{tabular}

        \caption{
          \it The optimal lead distance (in feet) by count and prior disengagements under the one-player MDP, assuming an average runner on first base facing an average battery (pitcher and catcher) with no outs and third base empty. Estimates include standard errors estimated via bootstrapping.
        }
        \label{tab:policy-by-count-mdp}
      \end{table}

      For the most part, the trends in the optimal lead distance policy in the one-player MDP are similar to that of the two-player game. As in the two-player game, the optimal lead distance primarily depends on prior disengagements, not count. The primary difference between Tables \ref{tab:policy-by-count-zsg} and \ref{tab:policy-by-count-mdp} is that in the one-player MDP the optimal lead is shorter with two prior disengagements. Another small difference is that in the one-player MDP the optimal lead distance increases in three-ball counts because pitchers are less likely to attempt pickoffs (Table \ref{tab:model-summary}).

      The increase in recommended lead distance following an unsuccessful pickoff attempt follows a simple pattern. Table \ref{tab:policy-increase-mdp} reports the difference in optimal lead distance between each state and its counterpart with one fewer disengagement. This shows how much a runner should increase their lead after an unsuccessful pickoff attempt, depending on the count, outs, and disengagement number. Weighted by the frequency with each unsuccessful pickoff attempts occur in each state, the average lead increase after the first disengagement is 1.65 feet and varies little by count and outs. After the second disengagement, the average lead increase is 2.31 feet and varies little by count and outs.
      
      \begin{table}
        \small
        \centering
        \begin{tabular}{c|ccc|ccc}
& \multicolumn{3}{c}{Increase after 1$^\text{st}$ Disengagement} & \multicolumn{3}{c}{Increase after 2$^\text{nd}$ Disengagement}\\
Count & 0 Outs & 1 Out & 2 Outs & 0 Outs & 1 Out & 2 Outs \\
  \hline
0-0 & 1.6 {\color{gray}$\pm$ 0.1} & 1.7 {\color{gray}$\pm$ 0.1} & 1.6 {\color{gray}$\pm$ 0.1} & 2.6 {\color{gray}$\pm$ 0.2} & 2.5 {\color{gray}$\pm$ 0.2} & 2.4 {\color{gray}$\pm$ 0.2} \\ 
  1-0 & 1.7 {\color{gray}$\pm$ 0.1} & 1.8 {\color{gray}$\pm$ 0.1} & 1.6 {\color{gray}$\pm$ 0.1} & 2.4 {\color{gray}$\pm$ 0.2} & 2.4 {\color{gray}$\pm$ 0.2} & 2.4 {\color{gray}$\pm$ 0.2} \\ 
  2-0 & 1.6 {\color{gray}$\pm$ 0.2} & 1.8 {\color{gray}$\pm$ 0.1} & 1.5 {\color{gray}$\pm$ 0.2} & 2.4 {\color{gray}$\pm$ 0.2} & 2.4 {\color{gray}$\pm$ 0.3} & 2.1 {\color{gray}$\pm$ 0.4} \\ 
  3-0 & 1.6 {\color{gray}$\pm$ 0.1} & 1.8 {\color{gray}$\pm$ 0.1} & 1.5 {\color{gray}$\pm$ 0.2} & 2.0 {\color{gray}$\pm$ 0.3} & 2.0 {\color{gray}$\pm$ 0.3} & 1.7 {\color{gray}$\pm$ 0.6} \\ 
   \hline
0-1 & 1.7 {\color{gray}$\pm$ 0.1} & 1.7 {\color{gray}$\pm$ 0.1} & 1.6 {\color{gray}$\pm$ 0.1} & 2.6 {\color{gray}$\pm$ 0.2} & 2.4 {\color{gray}$\pm$ 0.2} & 2.4 {\color{gray}$\pm$ 0.2} \\ 
  1-1 & 1.6 {\color{gray}$\pm$ 0.1} & 1.7 {\color{gray}$\pm$ 0.2} & 1.6 {\color{gray}$\pm$ 0.1} & 2.5 {\color{gray}$\pm$ 0.2} & 2.5 {\color{gray}$\pm$ 0.2} & 2.4 {\color{gray}$\pm$ 0.2} \\ 
  2-1 & 1.7 {\color{gray}$\pm$ 0.1} & 1.7 {\color{gray}$\pm$ 0.2} & 1.5 {\color{gray}$\pm$ 0.1} & 2.5 {\color{gray}$\pm$ 0.2} & 2.5 {\color{gray}$\pm$ 0.2} & 2.4 {\color{gray}$\pm$ 0.2} \\ 
  3-1 & 1.6 {\color{gray}$\pm$ 0.2} & 1.7 {\color{gray}$\pm$ 0.1} & 1.5 {\color{gray}$\pm$ 0.2} & 2.1 {\color{gray}$\pm$ 0.3} & 2.1 {\color{gray}$\pm$ 0.3} & 1.9 {\color{gray}$\pm$ 0.4} \\ 
   \hline
0-2 & 1.6 {\color{gray}$\pm$ 0.1} & 1.7 {\color{gray}$\pm$ 0.1} & 1.6 {\color{gray}$\pm$ 0.1} & 2.5 {\color{gray}$\pm$ 0.2} & 2.5 {\color{gray}$\pm$ 0.2} & 2.3 {\color{gray}$\pm$ 0.2} \\ 
  1-2 & 1.7 {\color{gray}$\pm$ 0.2} & 1.6 {\color{gray}$\pm$ 0.1} & 1.7 {\color{gray}$\pm$ 0.1} & 2.5 {\color{gray}$\pm$ 0.2} & 2.3 {\color{gray}$\pm$ 0.2} & 2.3 {\color{gray}$\pm$ 0.2} \\ 
  2-2 & 1.6 {\color{gray}$\pm$ 0.2} & 1.6 {\color{gray}$\pm$ 0.2} & 1.5 {\color{gray}$\pm$ 0.2} & 2.3 {\color{gray}$\pm$ 0.2} & 2.3 {\color{gray}$\pm$ 0.3} & 2.2 {\color{gray}$\pm$ 0.3} \\ 
  3-2 & 1.5 {\color{gray}$\pm$ 0.2} & 1.6 {\color{gray}$\pm$ 0.2} &  & 2.0 {\color{gray}$\pm$ 0.3} & 1.9 {\color{gray}$\pm$ 0.3} &  \\ 
\end{tabular}

        \caption{
          \it The increase in optimal lead distance (in feet) after an unsuccessful pickoff attempt under the one-player MDP, assuming an average runner on first base facing an average battery (pitcher and catcher) with third base empty. Estimates include standard errors estimated via bootstrapping.
        }
        \label{tab:policy-increase-mdp}
      \end{table}
  
      The qualitative findings from this section generalize approximately when extending beyond the average runner and battery. Table \ref{tab:policy-increase-by-skill-mdp} shows the optimal runner policy by number of disengagements for different combinations of runner skill and battery skill, along with the average increase in lead after each disengagement. Against batteries more skilled at controlling the run game, the optimal lead distance is shorter. For more skilled runners, the optimal lead distance is longer. Counterintuitively, we recommend greater lead increases for less skilled runners against more skilled batteries. This is because the recommended lead distance with zero disengagements is shorter for less skilled runners against more skilled batteries. When the runner has a greater skill advantage, the recommended lead starts longer. Depending on the skills of the players involved, the average recommended lead distance increase after the first disengagement ranges from 1.3 to 2.0 feet. After the second disengagement, the recommended increase ranges from 1.8 to 3.3 feet.
    
      \begin{table}
        \resizebox{\linewidth}{!}{
          \begin{tabular}{cc|ccc|cc}
\multicolumn{2}{c|}{Skill Percentile} & \multicolumn{3}{c|}{Lead by Disengagements (0 outs, 0-0 count)} & \multicolumn{2}{c}{Mean Increase (all outs/counts)}\\
Runner & Battery & 0 & 1 & 2 & after 1$^\text{{st}}$ Dis. & after 2$^\text{{nd}}$ Dis. \\
  \hline
90$^\text{th}$ & 10$^\text{th}$ & 12.1 {\color{gray}$\pm$ 0.4} & 13.4 {\color{gray}$\pm$ 0.4} & 15.3 {\color{gray}$\pm$ 0.5} & 1.3 {\color{gray}$\pm$ 0.1} & 1.8 {\color{gray}$\pm$ 0.3} \\ 
  90$^\text{th}$ & 50$^\text{th}$ & 10.7 {\color{gray}$\pm$ 0.3} & 12.2 {\color{gray}$\pm$ 0.3} & 14.2 {\color{gray}$\pm$ 0.4} & 1.5 {\color{gray}$\pm$ 0.1} & 1.9 {\color{gray}$\pm$ 0.3} \\ 
  90$^\text{th}$ & 90$^\text{th}$ & 8.8 {\color{gray}$\pm$ 0.6} & 10.4 {\color{gray}$\pm$ 0.5} & 12.7 {\color{gray}$\pm$ 0.5} & 1.6 {\color{gray}$\pm$ 0.1} & 2.1 {\color{gray}$\pm$ 0.2} \\ 
   \hline
50$^\text{th}$ & 10$^\text{th}$ & 11.1 {\color{gray}$\pm$ 0.4} & 12.6 {\color{gray}$\pm$ 0.4} & 14.7 {\color{gray}$\pm$ 0.4} & 1.5 {\color{gray}$\pm$ 0.1} & 2.0 {\color{gray}$\pm$ 0.3} \\ 
  50$^\text{th}$ & 50$^\text{th}$ & 9.6 {\color{gray}$\pm$ 0.3} & 11.3 {\color{gray}$\pm$ 0.3} & 13.7 {\color{gray}$\pm$ 0.3} & 1.6 {\color{gray}$\pm$ 0.1} & 2.2 {\color{gray}$\pm$ 0.2} \\ 
  50$^\text{th}$ & 90$^\text{th}$ & 7.6 {\color{gray}$\pm$ 0.6} & 9.5 {\color{gray}$\pm$ 0.5} & 12.2 {\color{gray}$\pm$ 0.6} & 1.8 {\color{gray}$\pm$ 0.1} & 2.5 {\color{gray}$\pm$ 0.2} \\ 
   \hline
10$^\text{th}$ & 10$^\text{th}$ & 9.7 {\color{gray}$\pm$ 0.4} & 11.3 {\color{gray}$\pm$ 0.4} & 14.1 {\color{gray}$\pm$ 0.4} & 1.6 {\color{gray}$\pm$ 0.1} & 2.5 {\color{gray}$\pm$ 0.2} \\ 
  10$^\text{th}$ & 50$^\text{th}$ & 8.3 {\color{gray}$\pm$ 0.3} & 10.2 {\color{gray}$\pm$ 0.3} & 13.3 {\color{gray}$\pm$ 0.3} & 1.8 {\color{gray}$\pm$ 0.1} & 2.9 {\color{gray}$\pm$ 0.2} \\ 
  10$^\text{th}$ & 90$^\text{th}$ & 6.5 {\color{gray}$\pm$ 0.6} & 8.6 {\color{gray}$\pm$ 0.6} & 12.0 {\color{gray}$\pm$ 0.6} & 2.0 {\color{gray}$\pm$ 0.2} & 3.3 {\color{gray}$\pm$ 0.3} \\ 
\end{tabular}

        }
        \caption{
          \it Average increase in optimal lead distance (in feet) by runner and battery skill percentile under the one-player MDP. The reported increase reflects the average increase across game states in which unsuccessful pickoff attempts occur. The hypothetical 90$^\text{th}$-percentile battery represents the 90$^\text{th}$ percentile at attempting pickoffs, successfully executing pickoffs, suppressing stolen base attempts, and catching would-be base stealers.
        }
        \label{tab:policy-increase-by-skill-mdp}
      \end{table}

  \section{Discussion and actionable insights}
  \label{sec:discussion}

    In comparing the optimal lead distance policy under the two-player game (Table \ref{tab:policy-by-count-zsg}) and the one-player MDP (Table \ref{tab:policy-by-count-mdp}), we observe substantial differences between them. For most states, the former lead distance is longer than the latter. The optimal policy for the runner under the two-player game is to extend their lead until the pitcher is indifferent between attempting a pickoff or not. By contrast, the one-player optimal policy corresponds to a predicted pickoff attempt probability below 10\%. Because the pitcher sometimes attempts low-success pickoffs, the runner does not need to extend their lead further.

    In the two-player game we recommend long leads (beyond 15 feet in some cases) after two prior disengagements. These distances are near the edge of the range of leads we observe but not outside the range. To develop intuition for this policy, consider the optimal lead distance of 14.4 feet in a 3-2 count after two prior disengagements with no outs (Table \ref{tab:policy-by-count-zsg}). If an average pitcher attempts a pickoff in this scenario, we estimate a 17.0\% probability that the runner is out. Otherwise, the runner advances to second base. By contrast, if the pitcher throws a pitch and if the runner attempts a stolen base, we estimate an 82.4\% probability that the stolen base is successful. In other words, the pitcher is indifferent to throwing a pitch or attempting a pickoff because the probability of an unsuccessful pickoff attempt is similar to the probability of a successful stolen base (conditioned on a steal attempt), and the outcome is the same either way.

    Because a third failed pickoff attempt results directly in free baserunner advancement, the pitcher is strongly disincentivized from attempting a pickoff after two prior disengagements. In the two-player game, this is reflected in a drastic increase in the optimal lead distance after the second disengagement (over 4 feet in most cases, Table \ref{tab:policy-by-count-zsg}). In the one-player MDP, we observe a more steady increase in optimal lead distance as a function of prior disengagements because observed pitcher pickoff attempt behavior steadily decreases as a function of prior disengagements (see Figure \ref{fig:prob-pickoff}).

    The value of the start-of-inning state under the two-player stochastic game is 0.4903 runs, meaning that the expected number of runs scored per inning is 0.4903. Stripping away pitcher agency under the one-player MDP increases run expectancy only slightly, to 0.4905 runs. If we strip away agency for the runner and fix the policy to be the observed average lead distance by game state in 2023, then the start-of-inning run expectation drops to 0.4862 runs. Over a full season of 162 games of nine innings each, the difference between observed and optimal behavior amounts to 5.7 more runs for the offense, with a bootstrapped standard error of 1.2 runs. To add 5.7 expected runs in free agency would cost roughly \$6 million \parencite{clemens_what_2021}, which speaks to a significant (but not outlandish) impact. This value may have increased since the last study we found to have estimated it.

    For runners and coaches seeking a strategy for the lead distance problem to implement today, we recommend the optimal policy from the one-player MDP rather than the two-player game. Pitchers are unlikely to immediately adapt their behavior to the optimal policy, and runners can exploit this. Because the pickoff attempt action is observable, runners can adjust their lead distance policy if they observe pitchers making different pickoff decisions over time.

    A simple approximation to the optimal policy is the \textit{Two-Foot Rule}: Extend your lead by two feet after each pickoff attempt. Because runners already have significant cognitive burden in professional baseball games, it is essential to boil down our results into a memorable and easily communicated rule of thumb. The Two-Foot Rule does not require consulting a chart to know what to do. It is a real-life actionable insight that would increase runner lead distances and increase run scoring.
   
  \section{Conclusion}
  \label{sec:conclusion}

    Under the new pickoff rules, baseball players and coaches understand intuitively that runners want to increase their lead distance after each successive disengagement by the pitcher---but by how much? In this paper, we formalize the cat-and-mouse game between pitcher and runner as a two-player sequential game and estimate the optimal lead distance policy for a runner on first base with other bases empty, given the number of outs, the ball-strike count, and the number of prior disengagements by the pitcher. We model the probabilities of runner outcomes (pickoff attempt/success, stolen base attempt/success) as functions of lead distance, context, and player skill. This allows us to estimate the optimal lead distance policy for any combination of pitcher skill and runner skill.

    We consider a two-player sequential game in which the runner first decides their lead distance and then the pitcher (upon observing the runner's lead distance) then decides whether to attempt a pickoff or throw a pitch. We show that the game-theoretic optimal policy for the runner is a maximin strategy (with respect to rest-of-inning run expectancy) so that the pitcher is indifferent to attempting a pickoff or throwing a pitch. We estimate model parameters using empirical data from the 2023 MLB season and find that under the game-theoretic equilibrium, runners would be more aggressive with their lead distances than the current behavior we observe. While baseball analytics sometimes garners negative attention for diminishing the aesthetic appeal of the game \parencite{sheinin_analytics_2025}, this is an example for which the analytically recommended strategy could be more aesthetically pleasing. In other sports, stochastic and dynamic models have shown that decision-makers have been traditionally too conservative \parencite{romer_firms_2006, beaudoin_strategies_2010, yam_what_2019}, but this is a new result for leadoff and pickoff decisions in baseball.

    We go a step further and consider the opportunity for runners to exploit the observed suboptimal behavior by pitchers. We consider a one-player MDP in which the pitcher's agency in the pickoff decision is replaced by a stochastic process, modeled on empirical data as a function of the runner's lead distance.  This framework leads to the most actionable recommendation for runners under the current conditions of pitcher behavior. We estimate that the average team could expect to increase their expected offensive output by approximately 5.7 runs over a full season by implementing the optimal lead distance policy. Based on these results, we propose the Two-Foot Rule: Extend your lead by two feet after each pickoff attempt. This easily implemented recommendation is a good approximation to the optimal lead distance policy across a wide array of contexts (outs and ball-strike count) and pitcher/runner skills.

    We carefully considered the lead distance problem and chose several assumptions to make it tractable. Most notably, we assumed that the runner's decision to attempt a stolen base is stochastic and independent of lead distance. We did not have appropriate data to model the true dynamic of this decision because it depends on latent variables like how the runner reads the body language of the pitcher. We also assumed that observed lead distance fully describes the runner's behavior and that they cannot disguise their aggressiveness by leaning in one direction. MLB has begun collecting 18-point pose tracking data on runners and fielders \parencite{jedlovec_introducing_2020}, which would enable more granular modeling of the runner-pitcher game if made available. With these data, future work might reveal a mixed strategy to be optimal for the pitcher's pickoff decision because the runner could hide their action. We assumed that lead distance only impacts pickoff and stolen base outcomes, not the outcomes of batted balls. Baseball domain knowledge suggests this is a reasonable assumption because while a pitch is being thrown, the runner moves to a secondary lead distance which is more important than primary lead distance for determining advancement outcomes. We assumed that, other than pickoff attempts, no other disengagements occur. In reality, the pitcher may derive value from stepping off the pitching rubber if they are short of breath or need to resolve a miscommunication with the catcher. These extensions are opportunities for future work. Another future direction would be to model an entire roster, which would involve a nonstationary two-player sequential game and introduce enormous technical challenges in analyzing and solving such a class of games.

    Ultimately, this paper addressed a timely problem of substantive interest in the baseball industry, achieving meaningful results. By leveraging compelling tracking data to model individual behavior and distill practical insights, our approach is related to a broader class of problems. In addition to sports leagues considering rule changes, this framework is relevant to any setting in which a regulator controls the rules under which competitors operate. Where high-fidelity data are available regarding competitor behavior, a similar approach can help regulators anticipate the consequences of rational competitor behavior.

\section*{Acknowledgments}
  Major League Baseball trademarks and copyrights are used with permission of MLB Advanced Media, L.P. All rights reserved. The authors thank Mallesh Pai of Rice University for constructive discussions during conceptualization of the problem and early implementation of the solution.

  \section*{Code and data availability}

    Code related to this paper is available at \url{https://github.com/jfhahn2/pickoff-game-theory}. The raw play-by-play lead distance data are proprietary but all other data and all code are open-source. All other data are publicly available and may be downloaded using the R package sabRmetrics \parencite{powers_sabrmetrics_2025}.

    \printbibliography

 \section*{Technical results}
    We present a formal description of the sequential game described in Section 2 and the related optimality characterizations.
    
    Recall from Section 2, the action space of the runner and pitcher are denoted by $A_{\text{R}}$ and $A_\text{P}$, respectively. The sets $\mathbb{F}_\text{R}, \mathbb{K}_\text{P}$, and $\mathbb{F}_\text{P}$ defined below will be used to rigorously define a policy. 
        \begingroup
        \allowdisplaybreaks
        \begin{align*}
            \mathbb{F}_\text{R} &\defeq \left\{f: {\cal S} \to \cup_{s \in {\cal S}}A_{\text{R}}(s) \big| f(s) \in A_{\text{R}}(s)\right\}\\
            \mathbb{K}_\text{P} &\defeq \left\{(s,a_\text{R})\big| s \in {\cal S}, a_\text{R} \in A_\text{R}(s)\right\}\\
            \mathbb{F}_\text{P} &\defeq \left\{f: \mathbb{K}_\text{P} \to \cup_{(s,a_\text{R}) \in \mathbb{K}_\text{P}} A_\text{P}(s,a_\text{R}) \big| f(s,a_\text{R}) \in A_\text{P}(s,a_\text{R})\right\}.
        \end{align*}
        \endgroup
    
       \noindent The set $\mathbb{F}_\text{R}$ is the set of all deterministic decision rules for the runner. Recall that a deterministic decision rule for the pitcher takes as input a state and a runner action, and maps it to an allowable action for the pitcher. The set $\mathbb{K}_\text{P}$ includes the set of all possible state and runner action pairs. The set $\mathbb{F}_\text{P}$ then captures the set of all deterministic decision rules for the pitcher. A deterministic, Markovian policy for the runner is a sequence $\pi = (\pi_0,\pi_1,\ldots,)$ such that $\pi_t \in \mathbb{F}_\text{R}$ for all integers $t \geq 0$. A deterministic, stationary policy for the runner is a sequence $\pi = (\pi_0,\pi_1,\ldots,)$ such that $\pi_t = \bar\pi \in \mathbb{F}_\text{R}$ for all integers $t \geq 0$. That is a stationary policy employs the same decision rule in all time-stages. Similarly, a deterministic, Markovian policy for the pitcher is a sequence $\tilde\pi = (\tilde\pi_0,\tilde\pi_1,\ldots,)$ such that $\tilde\pi_t \in \mathbb{F}_\text{P}$ for all integers $t \geq 0$. Finally, a deterministic, stationary policy for the pitcher is a sequence $\tilde\pi = (\tilde\pi_0,\tilde\pi_1,\ldots,)$ such that $\tilde\pi_t = \hat\pi \in \mathbb{F}_\text{P}$ for all integers $t \geq 0$. Consistent with the notation used in Section 2, we denote the set of deterministic, stationary policies for the runner and pitcher by $\Pi_\text{R}$ and $\Pi_\text{P}$, respectively. Following this style of notation, we denote the set of deterministic, Markovian policies for the pitcher and runner by $\overline\Pi_\text{R}$ and $\overline\Pi_\text{P}$, respectively. Clearly, $\Pi_\text{R} \subseteq \overline\Pi_\text{R}$ and $\Pi_\text{P} \subseteq \overline\Pi_\text{P}$.
    
       The optimal value function is then defined as
       \begingroup
       \allowdisplaybreaks
       \begin{equation}
            \label{eqn:sg_val_extended}
            V^*(s) \defeq \max_{\pi_\text{R} \in \overline\Pi_\text{R}}\min_{\pi_\text{P} \in \overline\Pi_\text{P}}V^{\pi_\text{R},\pi_\text{P}}(s) \ \forall s \in S.
       \end{equation}
       \endgroup
       The policy space in the right-hand side of \eqref{eqn:sg_val_extended} is over the larger class of deterministic, Markovian policies rather than deterministic, stationary policies as in (6). With an abuse of notation, we overload $V^*$ to denote the optimal value function in (6) as well as in \eqref{eqn:sg_val_extended}. The reason for doing this is that under Assumptions \ref{assum:game_halts} and \ref{assum:transition_pby_continuous}, Proposition \ref{lem:stationary_pol_optimal} proves the existence of   deterministic, stationary policies for the runner and pitcher that attain the respective maximum and minimum in \eqref{eqn:sg_val_extended}. That is, the policy spaces in the right-hand side of \eqref{eqn:sg_val_extended} can be restricted to the smaller class of deterministic, stationary policies as in (6) without affecting the optimal value in \eqref{eqn:sg_val_extended}. 
        
       As alluded in Section (2), we assume that an inning eventually halts. This can be mathematically captured via the following condition.
       \begin{assumption}\label{assum:game_halts}
           There exists $m \in \mathbb{N}$ such that 
           \begingroup
           \allowdisplaybreaks
           \begin{align*}
                \rho &\defeq \sup_{\pi_\text{R} \in \overline\Pi_\text{R}, \pi_\text{P} \in \overline\Pi_\text{P}, s \in {\cal S}} \rho(\pi_\text{R},\pi_\text{P},s) < 1,\\
                \text{ where }       \rho(\pi_\text{R},\pi_\text{P},s) &\defeq \mathbb{P}^{\pi_\text{R}, \pi_\text{P}}[S_m \neq \Delta|S_0=s] \ \forall \pi_\text{R} \in \overline\Pi_\text{R}, \, \pi_\text{P} \in \overline\Pi_\text{P}, \, s \in {\cal S}.          
           \end{align*}
           \endgroup
       \end{assumption}
       \noindent Assumption \ref{assum:game_halts} states that regardless of the policies and the starting state, there is a positive probability (at least $1-\rho$) that an inning ends after no more than $m \in \mathbb{N}$ plays.
    
       Under Assumption \ref{assum:game_halts}, we prove that the expected number of runs scored by the batter, $V^{\pi_\text{R},\pi_\text{P}}$ is finite. Our proof is based on standard arguments used in the analyses of undiscounted problems \cite[Section 7.2]{bertsekas_dynamic_2012}. We present it here for completeness.  
    
       \begin{lemma}
           \label{lem:sg_pol_val_well_defined}
           Suppose Assumption \ref{assum:game_halts} holds. Then $V^{\pi_\text{R},\pi_\text{P}}(s)$ is well-defined for all $\pi_\text{R} \in \overline\Pi_\text{R}, \pi_\text{P} \in \overline\Pi_\text{P}$, and $s \in {\cal S}$. Also, $V^{\pi_\text{R},\pi_\text{P}}$ is uniformly bounded over $\overline\Pi_\text{R} \times \overline\Pi_\text{P} \times {\cal S}$.  
       \end{lemma}       
       \begin{proof}{Proof}
           Fix a $\pi_\text{R} \in \overline\Pi_\text{R}, \pi_\text{P} \in \overline\Pi_\text{P}$, and $s \in {\cal S}$. Then 
           \begingroup
           \allowdisplaybreaks
           \begin{align*}
               V^{\pi_\text{R},\pi_\text{P}}(s) &= \mathbb{E}^{\pi_\text{R},\pi_\text{P}}\left[\sum_{t=0}^\infty r(S_{t+1}|S_t) \Big| S_0 = s\right] = \lim_{t \to \infty} \mathbb{E}^{\pi_\text{R},\pi_\text{P}}\left[\sum_{j=0}^{t-1} r(S_{j+1}|S_j) \Big| S_0 = s\right],
           \end{align*}
           \endgroup
           where the second equality follows by the monotone convergence theorem since the rewards are nonnegative. Under Assumption \ref{assum:game_halts}, it follows by standard Markov chain arguments that
           \[
           \mathbb{P}^{\pi_\text{R}, \pi_\text{P}}[S_{km} \neq \Delta|S_0=s] \leq \rho^k \ \forall k \in \mathbb{N}.
           \]
           Hence, the expected total rewards earned between periods $km$ and $(k+1)m - 1$ is upper bounded by $m\rho^k \max_{s,s' \in {\cal S}}r(s'|s)$. Therefore
           \[
           V^{\pi_\text{R},\pi_\text{P}}(s) \leq \sum_{k=0}^\infty m\rho^k \max_{s,s' \in {\cal S}}r(s'|s) = \frac{m \max_{s,s' \in {\cal S}}r(s'|s)}{1-\rho},
           \]
           thereby establishing the existence of $V^{\pi_\text{R},\pi_\text{P}}(s)$. The uniform bound follows since the above upper bound on $V^{\pi_\text{R},\pi_\text{P}}(s)$ is independent of $\pi_\text{R},\pi_\text{P}$, and $s$.
       \end{proof}
       
        We present optimality characterizations in Theorem \ref{lem:vi_convergence} and Propositions \ref{lem:bellman_equations} to  \ref{lem:policy_iteration}. These results and their proofs are similar to the undiscounted problem in the single-player case \cite[Section 7.2]{bertsekas_dynamic_2012} and the zero-sum stochastic game with simultaneous moves \cite{patek_stochastic_1999}. We state and prove our results for completeness and to outline the minor differences. 
      \begin{definition}\label{defn:sg_operators}
        \normalfont
        Consider the sequential game defined in Section (2). Fix a $u \in \mathbb{R}^{\abs{\cal S}}$ with $u(\Delta) = 0$.
        \begin{enumerate}
            \item For any $\pi_\text{R} \in \Pi_\text{R}$ and $\pi_\text{P} \in \Pi_\text{P}$, the evaluation operator $T_{\pi_\text{R},\pi_{\text{P}}}$ on $\mathbb{R}^{\abs{\cal S}}$ is defined as
            \begingroup
            \allowdisplaybreaks
            \begin{align}
            \label{eqn:eval_operator}
            (T_{\pi_\text{R},\pi_{\text{P}}} u)(s) &\defeq \mathbb{E}_{s' \sim p(\cdot|s,\pi_\text{R}(s),\pi_\text{P}(s,\pi_\text{R}(s)))}\left[r(s'|s) + u(s')\right] \ \forall s \in {\cal S}.
            \end{align}
            \endgroup
            \item The optimality operator $T$ on $\mathbb{R}^{\abs{\cal S}}$ is defined as
            \begingroup
            \allowdisplaybreaks
            \begin{align}
            \label{eqn:bellman_operator}
            (T u)(s) &\defeq \max_{a_\text{R} \in A_\text{R}(s)}\min_{a_\text{P} \in A_\text{P}(s,a_\text{R})}\mathbb{E}_{s' \sim p(\cdot|s,a_\text{R},a_\text{P})}\left[r(s'|s) + u(s')\right] \ \forall s \in {\cal S}.
            \end{align}
            \endgroup
        \end{enumerate}       
       \end{definition}

          \begin{theorem}
           \label{lem:vi_convergence}
           Consider the sequential game defined in Section (2). Fix a $u_0 \in \mathbb{R}^{\abs{\cal S}}$ with $u_0(\Delta) = 0$. Suppose Assumption \ref{assum:game_halts} holds. Then
           \begin{enumerate}
               \item $\lim_{k \to \infty}u_k(s) = V^*(s)$ for all $s \in {\cal S}$, where $u_{k} \defeq Tu_{k-1} \ \forall k \in \mathbb{N}$.
               \item For any $\pi_\text{R} \in \Pi_\text{R}$ and $\pi_\text{P} \in \Pi_\text{P}$, we have $\lim_{k \to \infty}\tilde u_k(s) = V^{\pi_\text{R},\pi_\text{P}}(s)$ for all $s \in {\cal S}$, where $\tilde u_k \defeq T_{\pi_\text{R},\pi_\text{P}}\tilde u_{k-1} \ \forall k \in \mathbb{N}$ and $\tilde u_0 \defeq u_0$.
           \end{enumerate}
       \end{theorem}
       \begin{proof}{Proof}       
           Consider part $(a)$.  Let $\pi_\text{R} = (\mu_0,\mu_1,\ldots,)$ and $\pi_\text{P} = (\nu_0,\nu_1,\ldots)$ be any  deterministic, Markovian (may not be stationary) policies for the runner and pitcher, respectively. Let $m \in \mathbb{N}$ be as defined in Assumption \ref{assum:game_halts}. Fix $s_0 \in {\cal S}$. For any $K \in \mathbb{N}$, we have
           \begingroup
           \allowdisplaybreaks
           \begin{align}
           \nonumber
               V^{\pi_\text{R},\pi_\text{P}}(s_0) &= \mathbb{E}^{\pi_\text{R},\pi_\text{P}}\left[\sum_{t=0}^{mK-1} r(S_{t+1}|S_t) + \sum_{t=mK}^\infty r(S_{t+1}|S_t) \Big| S_0 = s_0\right]\\
               \nonumber
               &= \mathbb{E}^{\pi_\text{R},\pi_\text{P}}\left[\sum_{t=0}^{mK-1} r(S_{t+1}|S_t) + u_0(S_{mK-1}) \Big| S_0 = s_0\right] -\\           \label{eqn:vi_convergence_temp1}
               &\qquad \qquad \mathbb{E}^{\pi_\text{R},\pi_\text{P}}\left[ u_0(S_{mK-1})\Big| S_0 = s_0\right] + \mathbb{E}^{\pi_\text{R},\pi_\text{P}}\left[\sum_{t=mK}^\infty r(S_{t+1}|S_t) \Big| S_0 = s_0\right].
           \end{align}
           \endgroup
           Using similar arguments as in the proof of Lemma \ref{lem:sg_pol_val_well_defined}, it follows that
           \begin{equation}      \label{eqn:vi_convergence_temp2}    \mathbb{E}^{\pi_\text{R},\pi_\text{P}}\left[\sum_{t=mK}^\infty r(S_{t+1}|S_t) \big| S_0 = s_0\right] \leq \sum_{k=K}^\infty m \rho^k \max_{s,s' \in {\cal S}}r(s'|s) = m\rho^K\max_{s,s' \in {\cal S}}r(s'|s)/(1-\rho).
           \end{equation}
           Next, 
           \begin{align}
               \nonumber
               &\abs{\mathbb{E}^{\pi_\text{R},\pi_\text{P}}\left[ u_0(S_{mK-1}) \vert S_0 = s_0\right]} = \abs{\sum_{s' \in {\cal S}}\mathbb{P}^{\pi_\text{R},\pi_\text{P}}[S_{mK-1}=s']u_0(s')}\\
               \nonumber
               &\qquad \leq \sum_{s' \in {\cal S}-\Delta}\mathbb{P}^{\pi_\text{R},\pi_\text{P}}[S_{mK-1}=s'] \left(\max_{s' \in {\cal S}-\Delta}\abs{u_0(s')}\right) + \mathbb{P}^{\pi_\text{R},\pi_\text{P}}[S_{mK-1}=\Delta]u_0(\Delta)\\
               \label{eqn:vi_convergence_temp3}
               &\qquad \leq \rho^K \max_{s' \in {\cal S}-\Delta}\abs{u_0(s')},
           \end{align}
           where the last inequality is due to Assumption \ref{assum:game_halts}. Utilizing \eqref{eqn:vi_convergence_temp2} and \eqref{eqn:vi_convergence_temp3} in \eqref{eqn:vi_convergence_temp1} gives us
           \begin{align*}
           &\mathbb{E}^{\pi_\text{R},\pi_\text{P}}\left[\sum_{t=0}^{mK-1} r(S_{t+1}|S_t) + u_0(S_{mK-1}) \Big| S_0 = s_0\right] -\rho^K \max_{s' \in {\cal S}-\Delta}\abs{u_0(s')} \leq  V^{\pi_\text{R},\pi_\text{P}}(s_0) \leq\\
           &\quad \mathbb{E}^{\pi_\text{R},\pi_\text{P}}\left[\sum_{t=0}^{mK-1} r(S_{t+1}|S_t) + u_0(S_{mK-1}) \Big| S_0 = s_0\right] +\rho^K \max_{s' \in {\cal S}-\Delta}\abs{u_0(s')} + \\
           &\qquad \qquad \frac{m\rho^K\max_{s,s' \in {\cal S}}r(s'|s)}{(1-\rho)}.
           \end{align*}
           From the definition of the operator $T$ in \eqref{eqn:bellman_operator}, it follows that $\max_{\pi_{\text{R}}}\min_{\pi_\text{P}}\mathbb{E}^{\pi_\text{R},\pi_\text{P}}\big[\sum_{t=0}^{mK-1}  \\ r(S_{t+1}|S_t) + u_0(S_{mK-1}) \big| S_0 = s_0\big] = u_k$. Hence, taking the minimum over all pitcher policies followed by taking the maximum over all runner policies in the above chain and then passing to limit $K \to \infty$ gives us $\lim_{K \to \infty}u_{Km}(s_0) = V^*(s_0)$. From Assumption \ref{assum:game_halts}, it follows that $\abs{u_{Km+t} - u_{Km}} \leq m \rho^K \max_{s,s' \in {\cal S}}r(s'|s)$ for all $t=1,\ldots,m$. This along with $\lim_{K \to \infty}u_{Km}(s_0) = V^*(s_0)$ establishes our claim that $\lim_{k \to \infty}u_k(s_0) = V^*(s_0)$.
    
           The proof of part $(b)$ follows via similar arguments as that of part $(a)$.
       \end{proof}
    
    Theorem \ref{lem:vi_convergence}(a) is value iteration for our sequential game. The computational experiments in Section 4 use value iteration to solve our sequential game after a suitable discretization of the runner's action space. Apart from serving as an algorithmic device, Theorem \ref{lem:vi_convergence} also enables proving the optimality of deterministic, stationary policies. This result established in Proposition \ref{lem:stationary_pol_optimal} uses the claim in Proposition \ref{lem:bellman_equations}, which in turn relies on the below assumption.
    
       \begin{assumption}
           \label{assum:transition_pby_continuous}
           The transition probability $p(\cdot|s,a_\text{R},a_\text{P})$ satisfies the following conditions. Fix $s \in {\cal S}$. Consider sequences $\{a_\text{R}^k\} \in A_\text{R}(s)$ and $\{a_\text{P}^k\} \in A_\text{P}(s,a_\text{R}^k)$ such that $\lim_{k \to \infty}a_\text{R}^k = a_\text{R} \in A_\text{R}(s)$ and  $\lim_{k \to \infty}a_\text{P}^k = a_\text{P} \in A_\text{P}(s,a_\text{R})$. Then $\lim_{k \to \infty}p(s'|s,a_\text{R}^k,a_\text{P}^k) = p(s'|s,a_\text{R},a_\text{P})$ for all $s' \in {\cal S}$.
       \end{assumption}
       \noindent Assumption \ref{assum:transition_pby_continuous} is a continuity condition imposed on the transition probabilities when viewed as a function of the runner and pitcher actions. Propositions  \ref{lem:bellman_equations} to \ref{lem:policy_iteration} rely upon this assumption to prove their stated claims. Our construction of transition probabilities in Section 3 satisfy this assumption. This assumption trivially holds if the action space of the runner $A_\text{R}(s)$ is finite for all $s \in {\cal S}$.
    
       \begin{proposition}
           \label{lem:bellman_equations}
           Consider the sequential game defined in Section (2). Suppose Assumptions \ref{assum:game_halts} and \ref{assum:transition_pby_continuous} hold. Then
           \begin{enumerate}
               \item The optimal value function $V^*$ is the unique solution to $u = Tu$, i.e., 
           \begingroup
           \allowdisplaybreaks
           \begin{equation}
               \label{eqn:bellman_equation}
               V^*(s) = \max_{a_\text{R} \in A_\text{R}(s)}\min_{a_\text{P} \in A_\text{P}(s,a_\text{R})}\mathbb{E}_{s' \sim p(\cdot|s,a_\text{R},a_\text{P})}\left[r(s'|s) + V^*(s')\right] \ \forall s \in {\cal S}.
           \end{equation}
           \endgroup
           \item For any $\pi_\text{R} \in \Pi_\text{R}$ and $\pi_\text{P} \in \Pi_\text{P}$, the function $V^{\pi_\text{R},\pi_\text{P}}$ is the unique solution to $u = T_{\pi_\text{R},\pi_\text{P}}u$, i.e., 
           \begingroup
           \allowdisplaybreaks
           \begin{equation}
           \label{eqn:bellman_eval}
               V^{\pi_\text{R},\pi_\text{P}}(s) =   \mathbb{E}_{s' \sim p(\cdot|s,\pi_\text{R}(s),\pi_\text{P}(s,\pi_\text{R}(s)))}\left[r(s'|s) + V^{\pi_\text{R},\pi_\text{P}}(s')\right] \ \forall s \in {\cal S}.          
           \end{equation}
           \endgroup
           \end{enumerate}
       \end{proposition}
       \begin{proof}{Proof}
           We first establish part $(a)$. Consider arbitrary $u_0 \in \mathbb{R}^{\abs{\cal S}}$ with $u_0(\Delta) = 0$ and define $u_{k} = Tu_{k-1} \ \forall k \in \mathbb{N}$. By Theorem \ref{lem:vi_convergence}, $\lim_{k \to \infty}u_k(s) = V^*(s)$ for all $s \in {\cal S}$. Passing to limit $k \to \infty$ in $u_k = Tu_{k-1}$,  we have for all $s \in {\cal S}$ 
           \begingroup
           \allowdisplaybreaks
           \begin{align}
           \nonumber
           V^*(s) &= \lim_{k\to\infty}\max_{a_\text{R} \in A_\text{R}(s)}\min_{a_\text{P} \in A_\text{P}(s,a_\text{R})}\mathbb{E}_{s' \sim p(\cdot|s,a_\text{R},a_\text{P})}\left[r(s'|s) + u_{k-1}(s')\right]\\
           \label{eqn:bellman_equation_temp1}
           &= \lim_{k \to \infty}\max_{a_\text{R} \in A_\text{R}(s)}\underbrace{\min_{a_\text{P} \in A_\text{P}(s,a_\text{R})}\sum_{s' \in {\cal S}}\left[r(s'|s) + u_{k-1}(s')\right]p(s'|s,a_\text{R},a_\text{P})}_{{\bf (I)}}.
           \end{align}
           \endgroup
           We prove that $\bf (I)$ converges (as $k \to \infty$) to $\min_{a_\text{P} \in A_\text{P}(s,a_\text{R})}\sum_{s' \in {\cal S}}\big[r(s'|s) + V^*(s')\big]p(s'|s,a_\text{R},a_\text{P})$ uniformly over $A_\text{R}(s)$ for each $s \in {\cal S}$.
           
           Fix a $s \in {\cal S}$. If $s = (b,c,o,d)$ is such that $b \neq (1,0,0)$, then from (1) and (2),  $A_\text{R}(s) = \{\delta\}$ and $ A_\text{P}(a_\text{R},s) = \{\delta\}$. As a result, $p(\cdot|s,a_\text{R},a_\text{P}) = p(\cdot|s)$, and hence the convergence is independent of $a_\text{R}$ and $a_\text{P}$. Suppose the component $b$ equals $(1,0,0)$. Let $\{a_\text{R}^k\} \in A_\text{R}(s)$ be such that $\lim_{k \to \infty}a_\text{R}^k = a_\text{R} \in A_\text{R}(s)$. Then from  (2), $A_\text{P}(a_\text{R}^k,s) = \{0,1\}$ for all $k \in \mathbb{N}$. Therefore
           \begingroup
           \allowdisplaybreaks
           \begin{align*}
           &\lim_{k \to \infty}\min_{a_\text{P} \in A_\text{P}(s,a_\text{R}^k)}\sum_{s' \in {\cal S}}\big[r(s'|s) + u_{k-1}(s')\big]p(s'|s,a_\text{R}^k,a_\text{P}) \\
           &= \lim_{k \to \infty}\min\Bigg(\sum_{s' \in {\cal S}}\big[r(s'|s) + u_{k-1}(s')\big]p(s'|s,a_\text{R}^k,0),\sum_{s' \in {\cal S}}\big[r(s'|s) + u_{k-1}(s')\big]p(s'|s,a_\text{R}^k,1)\Bigg)\\
           &= \min\Bigg(\sum_{s' \in {\cal S}}\big[r(s'|s) + V^*(s')\big]p(s'|s,a_\text{R},0),\sum_{s' \in {\cal S}}\big[r(s'|s) + V^*(s')\big]p(s'|s,a_\text{R},1)\Bigg)\\
           &= \min_{a_\text{P} \in A_\text{P}(s,a_\text{R})}\sum_{s' \in {\cal S}}\big[r(s'|s) + V^*(s')\big]p(s'|s,a_\text{R},a_\text{P}),
           \end{align*}
           \endgroup
           where the second equality follows from Assumption \ref{assum:transition_pby_continuous} and the convergence of $u_{k-1}$ to $V^*$ established in Theorem \ref{lem:vi_convergence}. Thus, we have established the continuous convergence of the sequence of functions $\big\{\min_{a_\text{P} \in A_\text{P}(s,a_\text{R})}\sum_{s' \in {\cal S}}\big[r(s'|s) + u_{k-1}(s')\big]p(s'|s,a_\text{R},a_\text{P})\big\}_k$ over $A_\text{R}(s)$. Next, the set $A_\text{R}(s)$ in (1)  is compact. Hence, by Theorem 5 in Subsection X, Section 21 in  \cite{kuratowski_topology_2014}, $\min_{a_\text{P} \in A_\text{P}(s,a_\text{R})}\sum_{s' \in {\cal S}}\big[r(s'|s) + u_{k-1}(s')\big]p(s'|s,a_\text{R},a_\text{P})$ converges to $\min_{a_\text{P} \in A_\text{P}(s,a_\text{R})}\sum_{s' \in {\cal S}}\big[r(s'|s) + V^*(s')\big]p(s'|s,a_\text{R},a_\text{P})$ uniformly over $A_\text{R}(s)$.
    
           Utilizing this uniform convergence in \eqref{eqn:bellman_equation_temp1} gives us 
           \begin{align*}
           &\lim_{k \to \infty}\max_{a_\text{R} \in A_\text{R}(s)}\min_{a_\text{P} \in A_\text{P}(s,a_\text{R})}\sum_{s' \in {\cal S}}\left[r(s'|s) + u_{k-1}(s')\right]p(s'|s,a_\text{R},a_\text{P})\\
           &= \max_{a_\text{R} \in A_\text{R}(s)}\lim_{k \to \infty}\min_{a_\text{P} \in A_\text{P}(s,a_\text{R})}\sum_{s' \in {\cal S}}\left[r(s'|s) + u_{k-1}(s')\right]p(s'|s,a_\text{R},a_\text{P})\\
           &= \max_{a_\text{R} \in A_\text{R}(s)}\min_{a_\text{P} \in A_\text{P}(s,a_\text{R})}\mathbb{E}_{s' \sim p(\cdot|s,a_\text{R},a_\text{P})}\left[r(s'|s) + V^*(s')\right],
           \end{align*}
           thereby establishing \eqref{eqn:bellman_equation}. The uniqueness of $V^*$ follows from the value iteration convergence result established in Theorem \ref{lem:vi_convergence}.
    
           Regarding part $(b)$, for each $s \in {\cal S}$, set $A_\text{R}(s) = \{\pi_\text{R}(s)\}$ and $A_\text{P}(s,\pi_\text{R}(s)) = \{\pi_\text{P}(s,\pi_\text{R}(s))\}$. The claim then  follows from part $(a)$.
       \end{proof}
    
       \begin{proposition}
           \label{lem:stationary_pol_optimal}
           Consider the sequential game defined in Section 2. Suppose Assumptions \ref{assum:game_halts} and \ref{assum:transition_pby_continuous} hold. Then a stationary pair $(\pi_\text{R}^*,\pi_\text{P}^*) \in \Pi_\text{R} \times \Pi_\text{P}$ is an equilibrium policy if and only if $\pi_\text{R}^*(s)$ and $\pi_\text{P}^*(s,\pi_\text{R}^*(s))$ attain the maximum and minimum in \eqref{eqn:bellman_equation} for all $s \in {\cal S}$.
       \end{proposition}
       \begin{proof}{Proof}
           Suppose $\pi_\text{R}^*(s)$ and $\pi_\text{P}^*(s,\pi_\text{R}^*(s))$ attain the maximum and minimum in \eqref{eqn:bellman_equation} for all $s \in {\cal S}$. Then \eqref{eqn:bellman_equation} becomes
           \[
           V^*(s) = \mathbb{E}_{s' \sim p(\cdot|s,\pi_\text{R}^*(s),\pi_\text{P}^*(s,\pi_\text{R}^*(s)))}\left[r(s'|s) + V^*(s')\right] = (T_{\pi_\text{R}^*,\pi_\text{P}^*}V^*)(s)  \ \forall s \in {\cal S},
           \]
           where the second equality follows by the definition of $T_{\pi_\text{R}^*,\pi_\text{P}^*}$ in \eqref{eqn:eval_operator}. By Proposition  \ref{lem:bellman_equations}(b), we then get $V^*(s)=V^{\pi_\text{R}^*,\pi_\text{P}^*}(s) \ \forall s \in {\cal S}$. Hence, $(\pi_\text{R}^*,\pi_\text{P}^*)$ is an equilibrium policy. Conversely, suppose $(\pi_\text{R}^*,\pi_\text{P}^*)$ is an equilibrium policy. Then $V^*(s) = V^{\pi_\text{R}^*,\pi_\text{P}^*}(s) \ \forall s \in {\cal S}$. From parts $(a)$ and $(b)$ of Proposition   \ref{lem:bellman_equations}, it follows that $\pi_\text{R}^*(s)$ and $\pi_\text{P}^*(s,\pi_\text{R}^*(s))$ attain the maximum and minimum in \eqref{eqn:bellman_equation} for all $s \in {\cal S}$.
       \end{proof}
    
    The next result states and proves the policy iteration algorithm for our two-player zero-sum sequential game.
       \begin{proposition}
           \label{lem:policy_iteration}
           Consider the sequential game defined in Section 2 and suppose Assumptions \ref{assum:game_halts} and \ref{assum:transition_pby_continuous} hold. Starting from any $\pi_\text{R}^0 \in \Pi_\text{R}$, define the following updates for all $k \geq 0$.
           \begin{align*}
              V^{\pi_{\text{R}}^k}(s) &\defeq \min_{\pi_\text{P} \in \Pi_\text{P}} V^{\pi_\text{R}^k,\pi_\text{P}}(s) \, \forall s \in \mathcal{S}, \mbox{ and}\\
              \label{eqn:update-policy-two-player}
              \pi_\text{R}^{k+1}(s) &\in \argmax_{a_\text{R} \in A_\text{R}(s)} \left\{\min_{a_\text{p} \in A_\text{P}(s,a_\text{R})}\sum_{s' \in \mathcal{S}}[r(s'|s) + V^{\pi_{\text{R}}^k}(s')] p(s'|s,a_\text{R},a_\text{P})\right\} \,  \forall s \in \mathcal{S}.
            \end{align*}
            Then $V^{\pi_\text{R}^k} \uparrow V^*$. 
       \end{proposition}
       \begin{proof}{Proof}
           For any $\pi_\text{R} \in \Pi_\text{R}$, define the operator $T_{\pi_\text{R}}$ on $\mathbb{R}^{\abs{\cal S}}$ as 
           \begingroup
           \allowdisplaybreaks
           \begin{equation}
               \label{eqn:policy_iteration_temp1}
               (T_{\pi_\text{R}} u )(s) \defeq \min_{a_\text{P} \in A_\text{P}(s,\pi_\text{R}(s))}\sum_{s' \in S}\left[r(s'|s) + u(s')\right]p(s'|s,\pi_\text{R}(s),a_\text{P}) \ \forall s \in {\cal S}.
           \end{equation}
           \endgroup
           Since $\pi_\text{R}$ is stationary, the best response of the pitcher corresponding to $\pi_\text{R}$ reduces to solving a MDP. Denote the value of this MDP as $V^{\pi_\text{R}}$, i.e., $V^{\pi_\text{R}}(s) \defeq \min_{\pi_\text{P} \in \Pi_\text{P}} V^{\pi_\text{R},\pi_\text{P}}(s) \ \forall s \in {\cal S}$. Using similar arguments as in the proof of Theorem  \ref{lem:vi_convergence}, it follows that $\lim_{k \to \infty}T_{\pi_\text{R}}^k u = V^{\pi_\text{R}}$, where $T_{\pi_\text{R}}^k$ is the composition of $T_{\pi_\text{R}}$ taken with itself $k$ times.  Also, using arguments similar to the proof of Proposition  \ref{lem:bellman_equations}, it follows that $V^{\pi_\text{R}}$ is the unique fixed point of $T_{\pi_\text{R}}$.
    
           Note that the policy sequence $\pi_\text{R}^k$ generated by our algorithm satisfies  
           \begingroup
           \allowdisplaybreaks
           \begin{equation}
           \label{eqn:pi_iteration_temp2}
           T_{\pi_\text{R}^{k+1}}V^{\pi_\text{R}^k} = TV^{\pi_\text{R}^k} \ \forall k \in \mathbb{N} \cup \{0\}.
           \end{equation}
           \endgroup
           Next, we observe that the sequence of functions $V^{\pi_\text{R}^k}$ correspond to the optimal value function of the pitcher's best response MDP. Hence, $V^{\pi_\text{R}^k} = T_{\pi_\text{R}^k} V^{\pi_\text{R}^k}$ for all integers $k \geq 0$. Since $Tu \geq T_{\pi_\text{R}}u$ for all $\pi_\text{R} \in \Pi_\text{R}$ and $u \in \mathbb{R}^{\abs{\cal S}}$, it follows that $V^{\pi_\text{R}^k} = T_{\pi_\text{R}^k} V^{\pi_\text{R}^k} \leq T V^{\pi_\text{R}^k}$ for all integers $k \geq 0$. Utilizing this in \eqref{eqn:pi_iteration_temp2} gives us
           \[
               T_{\pi_\text{R}^{k+1}}V^{\pi_\text{R}^k} = TV^{\pi_\text{R}^k} \geq V^{\pi_\text{R}^k} \ \forall k \in \mathbb{N} \cup \{0\}.
           \]
          Repeatedly applying $T_{\pi_\text{R}^{k+1}}$ in the above chain of expression, we get for any $j \in \mathbb{N}$ that $T_{\pi_\text{R}^{k+1}}^j V^{\pi_\text{R}^k} \geq TV^{\pi_\text{R}^k} \geq V^{\pi_\text{R}^k} \ \forall k \in \mathbb{N} \cup \{0\}$. Passing to limit $j \to \infty$ and noting that $\lim_{j \to \infty}T_{\pi_\text{R}^{k+1}}^j V^{\pi_\text{R}^k} = V^{\pi_\text{R}^{k+1}}$ gives us
           \begingroup
           \allowdisplaybreaks
           \begin{equation}
           \label{eqn:policy_iteration_temp2}
               V^{\pi_\text{R}^{k+1}} \geq TV^{\pi_\text{R}^k} \geq V^{\pi_\text{R}^k} \ \forall k \in \mathbb{N} \cup \{0\}.
           \end{equation}
           \endgroup
           By Lemma  \ref{lem:sg_pol_val_well_defined}, it follows that $V^{\pi_\text{R}}$ is uniformly bounded over $\Pi_\text{R} \times {\cal S}$. Hence, $V^{\pi_\text{R}^k} \uparrow \bar V \in \mathbb{R}^{\abs{\cal S}}$. We prove that $\bar V = V^*$. Using arguments similar to the ones used in the proof of Proposition  \ref{lem:bellman_equations}, it follows that $T$ is a  continuous operator on $\mathbb{R}^{\abs{\cal S}}$, i.e., $T u^k \to T u$ for all $u^k \to u$. Hence, passing to limit $k \to \infty$ in \eqref{eqn:policy_iteration_temp2} gives us $T \bar V = \bar V$. From Proposition  \ref{lem:bellman_equations}(a), it follows that $\bar V = V^*$, thereby completing the proof.
       \end{proof}
    
       If the state- and action-spaces are all finite, the convergence in th policy iteration algorithm happens in a finite number of steps. In that case, the policy obtained in the last iterate, say, $\pi_\text{R}^*$, is optimal for the runner. The pitcher's optimal policy (best response), $\pi_\text{P}^*$, can then be obtained as 
       \[
       \pi_\text{P}^*(s,\pi_\text{R}^*(s)) \in \argmin_{a_\text{P} \in A_\text{P}(s,\pi_\text{R}^*(s))}\sum_{s' \in {\cal S}}[r(s'|s) + V^*(s')] p(s'|s,\pi_\text{R}^*(s),a_\text{P}) \ \forall s \in {\cal S}, \pi_\text{R}^*(s) \in A_\text{R}(s).
       \]
 
\end{document}